\documentclass[11pt,reqno]{amsart}
\usepackage{amsfonts,amsmath,amsthm,amssymb,bm}
\usepackage[usenames,dvipsnames]{color}
\usepackage{enumerate,nicefrac}
\usepackage{hyperref,multicol}
\usepackage{rotating}
\usepackage[labelsep=space]{caption}
\usepackage[all]{xy}
\usepackage{tikz}
\usetikzlibrary{decorations.markings}
\usepackage[labelformat=empty]{caption}
\usepackage{mathtools}

\DeclarePairedDelimiter\floor{\lfloor}{\rfloor}
\newtheorem{theorem}{Theorem}[section]
\newtheorem{lemma}[theorem]{Lemma}

\newtheorem{corollary}[theorem]{Corollary}

\theoremstyle{definition}

\usepackage[ansinew]{inputenc}
\usepackage{ifthen}
\usepackage{animate}

\newcommand{\be}{\begin{equation}}
\newcommand{\ee}{\end{equation}}

\title[Representations of integers]{On the representation of an integer in Ostrowski and recurrence 
numeration systems}

\author{Mohit Mittal}
\address{Department of Mathematics\\
Birla~Inst{i}tute~of~Technology~and~Science, Pilani 333\,031 \textsc{India}}
\email{mohit.mittal\symbol{64}pilani.bits-pilani.ac.in}

\author{Divyum Sharma}
\address{Department of Mathematics\\
Birla~Inst{i}tute~of~Technology~and~Science, Pilani 333\,031 \textsc{India}}
\email{divyum.sharma\symbol{64}pilani.bits-pilani.ac.in}

\begin{document}

\date{}

%\begin{document}
\begin{abstract}
 We provide an effective upper bound for positive integers with bounded Hamming weights with respect to both a linear recurrence numeration system and an Ostrowski-$\alpha$ numeration system, where $\alpha$ is a quadratic irrational. We prove a similar result 
 for the representation of an integer in
 two \textit{different} Ostrowski numeration systems.
\end{abstract}

\maketitle
%=====================================

\section{Introduction}

In 1971, Senge and Straus \cite{SeSt1,SeSt2} proved that there are only finitely many integers with bounded sum of digits with respect to each of the bases $q_1$ and $q_2$ if and only if $\log q_1/\log q_2$ is irrational. In 1980, Stewart \cite{St_80}
proved an effective version of this result. As remarked by Bugeaud, Cipu and Mignotte \cite{BCM_13}, these results are  illustrations of the principle that two miracles cannot happen simultaneously for large integers in two \textit{unrelated} numeration systems. Luca and Bravo \cite{Luca_Bravo_2016} studied the problem of explicitly finding all integers with few non-zero digits in their representations in both the binary as well as in the Zeckendorf numeration systems (defined below). More precisely, they proved that 
the largest power of $2$ which can be written as the sum of two Fibonacci numbers is 16.  (See also
 \cite{Ziegler_Chim_2018,VZ_21}.)
Ziegler \cite{Zi_19} generalized Stewart's result to numeration systems with base sequence given by simple non-degenerate linear recurrence sequences with  \textit{dominant roots}.
The second author \cite{DS24} obtained an analogue of Stewart's result for canonical number systems in imaginary quadratic number fields. We refer to \cite{BHLS,EGST,Lu_00,Mi88,Sc90,Sp23} and the references therein for further results on this theme.

In this paper, we explore this theme for Ostrowski numeration systems and linear recurrence numeration systems. 
In 1922, Ostrowski \cite{Os_1922} proved that the sequence of denominators of the convergents to the simple continued fraction expansion of an irrational number forms the basis for a numeration system. More precisely, he established the following result (\cite[Theorem 3.9.1]{AS}).

%-----------
%\begin{theorem}\label{Ostro}
\textit{Let $\alpha$ be an irrational real number having continued fraction expansion $[a_0;a_1,\ldots]$. Let $(q_i)_{i\geq 0}$ be the sequence of the denominators of the convergents to the continued fraction expansion. Then every non-negative integer $n$ can be expressed uniquely as
	\be\label{decomp}
	n=\sum\limits_{0\leq i\leq \ell} \epsilon_{i} q_i,
	\ee
where the $\epsilon_{i}$'s are integers satisfying
	\begin{enumerate}[(i)]
	\item $0\leq \epsilon_{0}<a_1$.
	\item $0\leq \epsilon_{i}\leq a_{i+1}$ for $i\geq 1$.
	\item For $i\geq 1$, if $\epsilon_{i}=a_{i+1}$, then $\epsilon_{i-1}=0$.
	\end{enumerate}
}%\end{theorem}
%-----------
 Observe that Condition $(iii)$ ensures that the relation $q_{i+1}=a_{i+1}q_{i}+q_{i-1}$ cannot be used to replace a linear combination of summands with another summand. Further conditions $(i)-(iii)$ are equivalent to
 \begin{equation}\label{greediness}
     \sum_{i=0}^{j}\epsilon_{i} q_{i}< q_{j+1},
 \end{equation} 
 for all $ 0 \leq j \leq \ell$ (see \cite[Theorem 3.9.1]{AS}). The expression given in \eqref{decomp} is called the \textit{Ostrowski $\alpha$-representation} of $n$. If $\alpha=\nicefrac{(1+\sqrt{5})}{2}$, then the sequence of denominators is the sequence of Fibonacci numbers and the
Ostrowski $\alpha$-numeration system is called the \textit{Zeckendorf numeration system} \cite{Ze_72}. We refer to \cite{sur} for a survey on the connections between Ostrowski numeration systems and combinatorics of words, and to \cite{dyn} for a study of ergodic and topological-dynamical properties of various dynamical systems associated to Ostrowski $\alpha$-representations. We also refer to \cite{AST_022}, \cite{Be_72}, \cite{CRT_81}, \cite{Gel_68}, \cite{Kim_99}, \cite{DSMa'am_019}  for investigations on the joint distribution of the sum-of-digits functions associated with Ostrowski and radix numeration systems.

Let $\alpha$ be an irrational real number. For a non-negative integer $n$, let $H_{\alpha}(n)$ denote the Hamming weight, i.e. the number of non-zero digits in the 
Ostrowski-$\alpha$ representation of $n$. For a real quadratic irrational number  $\alpha$ with simple continued fraction expansion $\alpha=[a_0;\ldots,a_{r-1},\overline{b_0,\ldots,b_{s-1}}]$, $r\geq 0, s\geq 1$, we let 
\begin{equation}\label{t_alpha}
    t_{\alpha}=\textrm{trace} \left(\prod_{0\leq j<s}\begin{pmatrix}
    b_j & 1\\
    1 & 0
\end{pmatrix}\right).
\end{equation} 

We give an upper bound for the size of integers with bounded Hamming weights with respect to two \textit{different} Ostrowski numeration systems.
More precisely,
\begin{theorem}\label{TwoOst}
    Let $\alpha,\beta$ be real quadratic irrational numbers with simple continued fraction expansions 
    \begin{align*}
        \alpha&=[a_{0,0};\ldots,a_{0,r_0-1},\overline{b_{0,0},\ldots,b_{0,s_0-1}}], r_0\geq 0, s_0\geq 1,\\
        \beta&=[a_{1,0};\ldots,a_{1,r_1-1},\overline{b_{1,0},\ldots,b_{1,s_1-1}}], r_1\geq 0, s_1\geq 1
    \end{align*}
    and $\mathbb{Q}(\sqrt{t_{\alpha}^2-4(-1)^{s_0}}) \neq \mathbb{Q}(\sqrt{t_{\beta}^2-4(-1)^{s_1}})$. 
    Then there exists an effectively computable number $C$ depending only on $\alpha$ and $\beta$ such that if $n$ is a positive integer with
    \[
    H_{\alpha}(n)+H_{\beta}(n)\leq M,
    \]
    then 
    \[
    \log n\leq (CM\log M)^{M-1}.
    \]
\end{theorem}
We next consider the representation of non-negative integers with respect to linear recurrence sequences (defined below). We refer to \cite{Petho_Tichy_89} for further details.

\textit{Definition.} Let $G=(G_i)_{i \geq 0}$ be a sequence of positive integers which satisfy the relation $$G_{i+d}=e_1G_{i+d-1}+\cdots+e_d G_{i}, \ i \geq 0$$ and $e_{j} \in \mathbb{N} \cup \left\{0\right\}$ for $1 \leq j \leq d $ and $G_{0}=1$. This sequence is called a \textit{linear recurrence sequence} of order $d$ if 
\begin{enumerate}[(i)]
\item for $d=1$, we have $e_{1}>1$ and,
\item for $d \geq 2$, we have $e_{1} \geq \cdots \geq e_{d} >0$ and $G_{j} > e_{1}(G_{j-1}+\cdots+G_{0})$ where $1 \leq j \leq d-1$.
\end{enumerate}
If $\gamma_{1}, \gamma_{2}, \cdots,\gamma_{d}$ denote the 
roots of the characteristic polynomial 
\[
f(x)=x^{d}-e_{1}x^{d-1}-\cdots-e_{d}
\] 
of the linear recurrence sequence with $|\gamma_{1}| \geq |\gamma_{2}| \geq \cdots \geq  |\gamma_{d}|$,
then by \cite[Lemma 4]{Petho_Tichy_89}, we have $ e_{1} \leq \gamma_{1}< e_{1}+1$ and $|\gamma_{2}|<1$.\\
A linear recurrence sequence $G=(G_i)_{i \geq 0}$ forms the basis for a numeration system. For an arbitrary positive integer $n$, let $L$ be the integer such that $G_{L} \leq n < G_{L+1}$. Let $n_{L}=n$, $f_{j}=\floor{\nicefrac{n_{j}}{G_{j}}}$, $n_{j-1} = n_{j}-f_{j} G_{j}$ for $j=L, \ldots, 1$ and $f_{0}=n_{0}$. Then every positive integer has a well-defined representation of the form 
$$n =\sum_{j=0}^{L} f_{j}G_{j}.$$ We shall call this the base-$G$ representation of $n$. For a non-negative integer $n$, let $H_{G}(n)$ denote the Hamming weight, i.e. the number of non-zero digits in the base-$G$ representation of $n$. We prove

\begin{theorem}\label{OneOst}
    Let $\alpha$ be a real quadratic irrational number with simple continued fraction expansion  $$\alpha =[a_{0,0};\ldots,a_{0,r_0-1},\overline{b_{0,0},\ldots,b_{0,s_0-1}}], r_0\geq 0, s_0\geq 1$$ and let $G=(G_i)_{i \geq 0}$ be a linear recurrence sequence of order $d$  satisfying $G_{i+d}=e_1G_{i+d-1}+\cdots+e_d G_{i}$ for all $i \geq 0$, with  $\gamma_{1}$ as the real dominant root of its characteristic polynomial. Further assume that $\mathbb{Q}(\sqrt{t_{\alpha}^2 - 4(-1)^{s_0})}) \neq \mathbb{Q}\left(\gamma_{1}\right)$. Then there exists an effectively computable number $C'$ depending only on $\alpha$ and $G$ such that if $n$ is a positive integer with
    \[
    H_{\alpha}(n)+H_{G}(n)\leq M,
    \]
    then 
    \[
    \log n\leq (C'M\log M)^{M-1}.
    \]
\end{theorem}
Let $b$ be a positive integer $\geq 2$. Then $G=(b^{i})_{i=0}^{\infty}$ is a linear recurrence sequence 
with $d=1$ and $e_1 = b$. In this case, the base-$G$ representation of $n$ is the radix representation of $n$ in base-$b$ and, we denote the Hamming weight of $n$ by $H_{b}(n)$. We have the following corollary:
\begin{corollary}\label{OneOst_base_b}
    Let $\alpha$ be a real quadratic irrational number and let $b$ be a positive integer $\geq 2$. Then there exists an effectively computable number $C''$ depending only on $\alpha$ and $b$ such that if $n$ is a positive integer with
    \[
    H_{\alpha}(n)+H_{b}(n)\leq M,
    \]
    then 
    \[
    \log n\leq (C''M\log M)^{M-1}.
    \]
\end{corollary} 
Theorems \ref{TwoOst} \& \ref{OneOst} will be derived from the following technical theorems, which are more amenable to a proof by induction.
\begin{theorem}\label{TwoOstTech}
Let $\alpha, \beta$ be real quadratic irrational numbers with simple continued fraction expansions 
    \begin{align*}
        \alpha&=[a_{0,0};\ldots,a_{0,r_0-1},\overline{b_{0,0},\ldots,b_{0,s_0-1}}], r_0\geq 0, s_0\geq 1,\\
        \beta&=[a_{1,0};\ldots,a_{1,r_1-1},\overline{b_{1,0},\ldots,b_{1,s_1-1}}], r_1\geq 0, s_1\geq 1
    \end{align*} and $\mathbb{Q}(\sqrt{t_{\alpha}^2-4(-1)^{s_0}}) \neq \mathbb{Q}(\sqrt{t_{\beta}^2-4(-1)^{s_1}})$.
    Let $(q_{\alpha,i})_{i\geq 0},(q_{\beta,i})_{i\geq 0}$ be the sequences of denominators of the convergents to the continued fraction expansions of $\alpha,\beta$, respectively.
Then there exists a constant $C_0$, depending only on \\
$s_0,s_1,a_{0,0},\ldots,a_{0,r_0-1},b_{0,0},\ldots,b_{0,s_0-1},a_{1,0},\ldots,a_{1,r_1-1},b_{1,0},\ldots,b_{1,s_1-1}$, with the following property:\\
For every $u\geq 4$, there exist integers $K, L, K_4,\ldots, K_u, L_4,\ldots, L_u$ with
    \begin{align*}
        &2= K_4\leq\cdots\leq K_u=K,\\
        &2= L_4\leq\cdots\leq L_u=L,\\
        &K+L=u
    \end{align*}
such that if
    \begin{align*}
         n&=\epsilon_1q_{\alpha,N_1}+\cdots+\epsilon_kq_{\alpha,N_k},\ N_1>\cdots>N_k,\\
         n&=d_1q_{\beta,M_1}+\cdots+\epsilon_{\ell}q_{\beta,M_{\ell}},\ M_1>\cdots>M_{\ell}
    \end{align*}
are the Ostrowski-$\alpha$ and Ostrowski-$\beta$ representations of a positive integer $n$, respectively, with $N_1\geq 2s_0+r_0$, $\epsilon_{i} \neq 0$ for $1 \leq i \leq k$ and $d_{i} \neq 0$ for $1 \leq i \leq \ell$, 
then for each $i=4,\ldots,u$, we have 
    \begin{align}
    \label{tech_min}    \min(n_1-n_{K_i},m_1-m_{L_i})&\leq (C_0\log n_1)^{i-3}\\
    \textrm{ and }\ 
    \label{tech_max}    \max(n_1-n_{K_i-1},m_1-m_{L_i-1})&\leq (C_0\log n_1)^{i-4},
    \end{align}
    where
    \begin{align*}
      n_j:&=\floor*{\frac{N_j-r_0}{s_0}} \textrm{ for } 1\leq j\leq k,\  n_j:=0,\  \textrm{ for }  j>k,  \\
      m_j:&=\floor*{\frac{M_j-r_1}{s_1}}\  \textrm{ for } 1\leq j\leq \ell,\  m_j:=0,\  \textrm{ for } j>\ell.
    \end{align*}
\end{theorem}

\begin{theorem}\label{OneOstTech}
 Let $\alpha$ be a real quadratic irrational number with simple continued fraction expansion 
 \begin{equation*}
     \alpha=[a_{0,0};\ldots,a_{0,r_0-1},\overline{b_{0,0},\ldots,b_{0,s_0-1}}], r_0\geq 0, s_0\geq 1.
 \end{equation*}
    Let $(q_{\alpha,i})_{i\geq 0}$ be the sequence of denominators of the convergents to the continued fraction expansion of $\alpha$.
    Let $G=(G_i)_{i \geq 0}$ be a linear recurrence sequence of order $d$ satisfying $G_{i+d}=e_1G_{i+d-1}+\cdots+e_d G_{i}$, $i \geq 0$ with  $\gamma_{1}$ as the real dominant root of its characteristic polynomial. Further assume that $\mathbb{Q}(\sqrt{t_{\alpha}^2-4(-1)^{s_0}}) \neq \mathbb{Q}\left(\gamma_{1}\right)$. 
Then there exists a constant $C_{0}'$, depending only on 
$s_0,a_{0,0},\ldots,a_{0,r_0-1},b_{0,0},\ldots,b_{0,s_0-1},$ $e_{1}, \ldots, e_{d}$ with the following property:\\
For every $u\geq 4$, there exist integers $K, L, K_4,\ldots, K_u, L_4,\ldots, L_u$ with
    \begin{align*}
        &2= K_4\leq\cdots\leq K_u=K,\\
        &2= L_4\leq\cdots\leq L_u=L,\\
        &K+L=u
    \end{align*}
such that if
    \begin{align*}        n&=\epsilon_1q_{\alpha,N_1}+\cdots+\epsilon_kq_{\alpha,N_k},\ N_1>\cdots>N_k,\\
         n&=f_{1}G_{R_1}+\cdots+f_{s}G_{R_{s}},\ R_1>\cdots>R_{s}
    \end{align*}
are the Ostrowski-$\alpha$ and base-$G$ representations of a positive integer $n$, respectively, with $N_1\geq 2s_0+r_0$, $\epsilon_{i} \neq 0$ for $1 \leq i \leq k$ and $f_{i} \neq 0$ for $1 \leq i \leq s$,
then for each $i=4,\ldots,u$, we have 
    \begin{align}
    \label{tech_min_2}    \min(n_1-n_{K_i},R_1-R_{L_i})&\leq (C_{0}'\log n_1)^{i-3}\\
    \textrm{ and }\ 
    \label{tech_max_2}    \max(n_1-n_{K_i-1},R_1-R_{L_i-1})&\leq (C_{0}'\log n_1)^{i-4},
    \end{align}
    where
    \begin{align*}
      n_j:&=\floor*{\frac{N_j-r_0}{s_0}} \textrm{ for } 1\leq j\leq k,\  n_j:=0,\  \textrm{ for }  j>k,  \\
      R_j:&=0,\  \textrm{ for } j>s.
    \end{align*}
 
\end{theorem}

\textit{Remark.} The problem of bounding the size of integers with fixed Hamming weights with respect to two linear recurrence numeration systems was studied in \cite{Zi_19} provided both the base sequences have a \textit{dominant characteristic root}.
We note that if $\alpha$ is a real quadratic irrational number, then the sequence of denominators of the convergents of its continued fraction expansion indeed satisfies a linear recurrence relation, but its characteristic polynomial may not have a dominating root. For instance, if $\alpha=\sqrt{11}$, then its continued fraction expansion is $[3; \overline{3, 6}]$ and $q_{i+4} = 20q_{i+2}- q_{i}$, $i \geq 0$. The characteristic polynomial is $x^{4}-20x^{2}+1$ and its roots are $\pm \sqrt{10+3\sqrt{11}}, \pm \sqrt{10-3\sqrt{11}}$. To circumscribe this issue, we use an argument of Lenstra and Shallit \cite{LeSh93} to split the sequence into finitely many subsequences, each satisfying the same binary recurrence relation (with different initial terms). See Section \ref{prelim}.
We then use an induction argument (following \cite{Zi_19}) and Baker-type estimates for linear forms in logarithms to prove our theorems. 

In Section \ref{prelim}, we record some preliminary results. In Sections  \ref{Sec_TwoOstTech} and \ref{Sec_OneOstTech}, we prove Theorems \ref{TwoOstTech} and \ref{OneOstTech}, respectively.
In Section \ref{Sec_OneTwoOst}, we derive Theorems \ref{TwoOst} and \ref{OneOst}. In section \ref{sec_6}, we prove Corollary \ref{OneOst_base_b}. 

\section{Preliminaries}\label{prelim}

Let $\alpha$ be a real quadratic irrational number. Let its simple continued fraction expansion be $\alpha=[a_0;\ldots,a_{r-1},\overline{b_0,\ldots,b_{s-1}}]$, $r\geq 0, s\geq 1$.
Let $(q_{\alpha,i})_{i\geq 0}$ be the sequence of denominators of the convergents to $\alpha$. Then by \cite[p. 352]{LeSh93}, $(q_{\alpha,i})_{i\geq 0}$ satisfies the linear recurrence
    \begin{equation}\label{rr_qn}
    q_{\alpha,i+2s}-t_{\alpha}q_{\alpha,i+s}+(-1)^sq_{\alpha,i}=0,\ i\geq r,  
    \end{equation}
    where $t_{\alpha}$ is as defined in \eqref{t_alpha}.\\
%Note that $t\geq 2$ if $\alpha\notin \frac{-1+\sqrt{5}}{2}+\mathbb{Z}$.
The next lemma gives a Binet-type formula for the convergent denominators.
\begin{lemma}\label{Binet-type}
   Let $\alpha$ be a real quadratic irrational number with simple continued fraction expansion  $[a_0;\ldots,a_{r-1},\overline{b_0,\ldots,b_{s-1}}]$, $r\geq 0, s\geq 1$. For $j=0,\ldots,s-1$, let $q^{(j)}_{\alpha,i}=q_{\alpha,si+j+r},\ i\geq 0$. 
   \begin{enumerate}[(i)]
       \item We have \begin{equation}\label{binet_rep}
       q^{(j)}_{\alpha,i}=c_{\alpha,1}^{(j)}\theta_{\alpha,1}^{i}-c_{\alpha,2}^{(j)}\theta_{\alpha,2}^i,\ i \geq 0,
       \end{equation}
       where 
    \begin{align*}
        \theta_{\alpha,1}&=\frac{t_{\alpha}+\sqrt{t_{\alpha}^2-4(-1)^s}}{2},\  \theta_{\alpha,2}=\frac{t_{\alpha}-\sqrt{t_{\alpha}^2-4(-1)^s}}{2},\\
        c^{(j)}_{\alpha,1}&=\frac{q^{(j)}_{\alpha,1}-\theta_{\alpha,2}q^{(j)}_{\alpha,0}}{\theta_{\alpha,1}-\theta_{\alpha,2}},\ c_{\alpha,2}^{(j)}=\frac{q^{(j)}_{\alpha,1}-\theta_{\alpha,1}q^{(j)}_{\alpha,0}}{\theta_{\alpha,1}-\theta_{\alpha,2}}.
    \end{align*}
       \item There exist positive constants $c_{\alpha, 3}$, $c_{\alpha,4}'$, $N_{0}(\alpha)$ depending only on $\alpha$ such that for every $j=0,\ldots,s-1$,
       \begin{equation}\label{qn_ub_lb}
            q^{(j)}_{\alpha,i} \leq c_{\alpha,3} \theta_{\alpha,1}^{i}\ \textrm{ for }\ i \geq 0\ \textrm{ and }\ c_{\alpha,4}' \theta_{\alpha,1}^{i} \leq q^{(j)}_{\alpha,i} \textrm{ for }\ i \geq N_{0}(\alpha).
       \end{equation}
   \end{enumerate}
\end{lemma}
\begin{proof}
From \eqref{rr_qn} and the definition of $q^{(j)}_{\alpha,i}$, we have
    \[
    q^{(j)}_{\alpha,i+2}-t_{\alpha}q^{(j)}_{\alpha,i+1}+(-1)^sq^{(j)}_{\alpha,i}=0,\ i\geq 0.
    \]
Now it is easy to see that $(i)$ follows.\\
%     For $j=0,\ldots,s-1$, we have
%     \[
% q^{(j)}_{\alpha,i}=q_{\alpha,si+j+r},\ i\geq 0.
%     \]
% Then
%     \[
%     q^{(j)}_{\alpha,i+2}-t_{\alpha}q^{(j)}_{\alpha,i+1}+(-1)^sq^{(j)}_{\alpha,i}=0,\ i\geq 0.
%     \]
% Hence
%     \begin{equation}
% \nonumber q^{(j)}_{\alpha,i}=c_{\alpha,1}^{(j)}\theta_{\alpha,1}^{i}-c_{\alpha,2}^{(j)}\theta_{\alpha,2}^i,
%     \end{equation}
% where 
%     \begin{align*}
%         \theta_{\alpha,1}&=\frac{t_{\alpha}+\sqrt{t_{\alpha}^2-4(-1)^s}}{2},\  \theta_{\alpha,2}=\frac{t_{\alpha}-\sqrt{t_{\alpha}^2-4(-1)^s}}{2},\\
%         c^{(j)}_{\alpha,1}&=\frac{q^{(j)}_{\alpha,1}-\theta_{\alpha,2}q^{(j)}_{\alpha,0}}{\theta_{\alpha,1}-\theta_{\alpha,2}},\ c_{\alpha,2}^{(j)}=\frac{q^{(j)}_{\alpha,1}-\theta_{\alpha,1}q^{(j)}_{\alpha,0}}{\theta_{\alpha,1}-\theta_{\alpha,2}}.
%     \end{align*}
\textit{Proof of (ii)}.\\
We have
    \[
    t_{\alpha}-2\leq \sqrt{t_{\alpha}^2-4(-1)^s}\leq t_{\alpha}+2,
    \]
and hence
    \[
    \max(t_{\alpha}-1,1)< \theta_{\alpha,1}< t_{\alpha}+1,\ -1<\theta_{\alpha,2}< 1.
    \]
    The inequality holds because if $t_{\alpha}=1$, then $s=1$ and $\theta_{\alpha,1}=\nicefrac{(1+\sqrt{5})}{2}$.
    (Note that $\theta_{\alpha,1},\theta_{\alpha,2}$ are irrational.)
    Further, $c^{(j)}_{\alpha,1}>0$ for all $j=0,\ldots,s-1$. We get 
    \begin{equation}
    \nonumber q^{(j)}_{\alpha,i}\leq c_{\alpha,3}\theta_{\alpha,1}^i , \  i \geq 0,
    \end{equation}
where $c_{\alpha,3}=\max\limits_{0\leq j\leq s-1}(c^{(j)}_{\alpha,1}+|c^{(j)}_{\alpha,2}|)$. Note that $c_{\alpha,3}$ depends only on $\alpha$. Since $|\nicefrac{\theta_{\alpha,2}}{\theta_{\alpha,1}}|<1$, there exists $N_0(\alpha)\in\mathbb{N}$ depending only on $\alpha$ such that for $i\geq N_0(\alpha)$, we have
    \[
    \left|\frac{\theta_{\alpha,2}}{\theta_{\alpha,1}}\right|^i<\frac{c^{(j)}_{\alpha,1}}{2|c^{(j)}_{\alpha,2}|}.
    \]
Then for $i\geq N_0(\alpha)$, we obtain
    \begin{equation}
         \nonumber q^{(j)}_{\alpha,i}\geq\theta_{\alpha,1}^i\left||c^{(j)}_{\alpha,1}|-|c^{(j)}_{\alpha,2}|\left|\frac{\theta_{\alpha,2}}{\theta_{\alpha,1}}\right|^i\right|>\frac{c^{(j)}_{\alpha,1}}{2}\theta_{\alpha,1}^i\geq c_{\alpha,4}'\theta_{\alpha,1}^i,
    \end{equation}
where $c_{\alpha,4}'=\min\limits_{0\leq j\leq s-1}\nicefrac{c^{(j)}_{\alpha,1}}{2}$.
\end{proof}

%In what follows, $C_1,C_2,\ldots$ denote constants depending upon 
%$b,s,a_0,\ldots,a_{r-1},b_0,\ldots,b_{s-1}$ only. 

The following lemma compares a given positive integer $N$ with a suitable power of $\theta_{\alpha,1}$, depending on the highest index in its Ostrowski-$\alpha$ representation.

\begin{lemma}\label{N_ub_lb}
   Let $\alpha$ be a real quadratic irrational number with simple continued fraction expansion  $[a_0;\ldots,a_{r-1},\overline{b_0,\ldots,b_{s-1}}]$, $r\geq 0, s\geq 1$. If \[N = \epsilon_1q_{\alpha,N_1}+\cdots+\epsilon_kq_{\alpha,N_k},\ N_1\geq\cdots\geq N_k,\] is the Ostrowski-$\alpha$ representation of $N$ with $\epsilon_i\neq 0$ for all $i=1,\ldots,k$, then there exist positive constants $c_{\alpha, 4}$, $c_{\alpha, 5}$ depending only on $\alpha$ such that
   \begin{equation}\label{n_upperbound}
     c_{\alpha,4}\theta_{\alpha,1}^{n_1} \leq N< c_{\alpha,5}\theta_{\alpha,1}^{n_1}.
   \end{equation}
\end{lemma}
% Note that the constant $c_{\alpha,3}$ arising in the proof of the above lemma is the same as that of Lemma \ref{Binet-type}.
\begin{proof}  From \eqref{greediness}, we have
 \[
         N  <q_{\alpha,N_1+1}.
  \]
 If $N_1<r$, then 
  \[
  N<q_{\alpha,r}.
  \]
If not, write $N_i=sn_i+j_i+r$, $0\leq j_i\leq s-1$. Then
 \[
         N=  \sum_{i=1}^k \epsilon_iq^{(j_i)}_{\alpha,n_i}<q_{\alpha,N_1+1}.
  \]
Using Lemma \ref{Binet-type}, we get
  \[
  N<c_{\alpha,3}\theta_{\alpha,1}^{n_{1}+1}.
  %\leq c_{\alpha,5}\theta_{\alpha,1}^{n_{1}},
  \]
Thus, in all cases, $N<c_{\alpha,5}\theta_{\alpha,1}^{n_{1}}$    where $c_{\alpha,5}=\max(c_{\alpha,3}\theta_{\alpha,1},q_{\alpha,r})$.
Further, there exist constants $c_{\alpha,4}'$ and $N_0$, depending only on $\alpha$, such that  for $n_1\geq N_0$, we have 
    \begin{align*}
        c_{\alpha,4}'\theta_{\alpha,1}^{n_1} \leq q^{(j_1)}_{\alpha,n_1}\leq  N=  \sum_{i=1}^k \epsilon_iq^{(j_i)}_{\alpha,n_i}.
    \end{align*}
   Now, for $n_{1} < N_{0}$, we have $$c_{\alpha, 4}'' \theta^{n_1}_{\alpha,1}\leq N,$$ where $c_{\alpha, 4}''=\min_{n_{1} < N_{0}} \frac{q^{(j_1)}_{\alpha, n_{1}}}{\theta^{n_1}_{\alpha,1}}$.
   Combining, we get
 $$ c_{\alpha, 4} \theta^{n_1}_{\alpha,1} \leq q^{(j_1)}_{\alpha,n_1} \leq N,$$ where $c_{\alpha,4}=\min(c_{\alpha, 4}', c_{\alpha, 4}'')$.
\end{proof}
%  Suppose
%     \[
%    N= \epsilon_1q_{\alpha,N_1}+\cdots+\epsilon_kq_{\alpha,N_k},\ N_1\geq\cdots\geq N_k,
%     \]
%     is the Ostrowski-$\alpha$ representation of $N$ with $\epsilon_i\neq 0$ for all $i=1,\ldots,k$. Write
%     $N_i=sn_i+j_i+r$, $0\leq j_i\leq s-1$. Then for $n_1\geq N_0$, we have 
%     \begin{equation}\label{n_upperbound}
% c_{\alpha,4}\theta_{\alpha,1}^{n_1}<q^{(j_1)}_{\alpha,n_1}\leq  N=  \sum_{i=1}^k \epsilon_iq^{(j_i)}_{\alpha,n_i}<q_{\alpha,N_1+1}<c_{\alpha,5}\theta_{\alpha,1}^{n_1},
%     \end{equation}
%     where $c_{\alpha,5}=\max(c_{\alpha,3}\theta_{\alpha,1},1)$. The upper bound follows from ($\ast$).
 We recall that two complex numbers $v_{1}$ and $v_{2}$ are said to be \textit{multiplicatively independent} if for integers $t_{1}$ and $t_{2}$, $v_{1}^{t_{1}} v_{2}^{t_2}=1$ implies $t_{1}=t_{2}=0$.
The next lemma provides a sufficient condition for $\theta_{\alpha,1}$ and $\theta_{\beta,1}$ to be multiplicatively independent.
\begin{lemma}\label{TwoOst_rootsLI}
   Let $\alpha, \beta$ be real quadratic irrational numbers with simple continued fraction expansions 
    \begin{align*}
        \alpha&=[a_{0,0};\ldots,a_{0,r_0-1},\overline{b_{0,0},\ldots,b_{0,s_0-1}}], r_0\geq 0, s_0\geq 1,\\
        \beta&=[a_{1,0};\ldots,a_{1,r_1-1},\overline{b_{1,0},\ldots,b_{1,s_1-1}}], r_1\geq 0, s_1\geq 1.
    \end{align*} If $\mathbb{Q}(\sqrt{t_{\alpha}^2-4(-1)^{s_0}}) \neq \mathbb{Q}(\sqrt{t_{\beta}^2-4(-1)^{s_1}})$, then $\theta_{\alpha,1}$ and $\theta_{\beta,1}$ are multiplicatively independent.
\end{lemma}
\begin{proof}
    If $\theta_{\alpha,1}$ and $\theta_{\beta,1}$ are multiplicatively dependent, then there exist integers $k_{1}$ and $k_{2}$ such that $\theta_{\alpha,1}^{k_1} \theta_{\beta,1}^{k_2} =1.$ If $k_{2}=0$, then $\theta_{\alpha,1}^{k_1} =1$, which is not possible as $\theta_{\alpha,1} > 1$.
    % . Since $\theta_{\alpha,1} > \max (t_{\alpha}-1, 1)$, therefore, either $\theta_{\alpha,1}^{k_1}>1$ or $\theta_{\alpha,1}^{-k_1}>1$, hence $ k_{2} \neq 0$.  
    Similarly, $k_{1} \neq 0$. Note that $k_{1}$ and $k_{2}$ cannot both be positive or both be negative. Without loss of generality, assume $k_{1} < 0$ and $k_{2}>0$.
    Now, consider the embedding of 
$ \mathbb{Q}(\sqrt{t_{\alpha}^2-4(-1)^{s_0}}, \sqrt{t_{\beta}^2-4(-1)^{s_1}})$ which takes
\begin{align*}
    \sqrt{t_{\alpha}^2-4(-1)^{s_0}} \ &\textrm{ to } \  \sqrt{t_{\alpha}^2-4(-1)^{s_0}}\\
    \noindent \textrm{ and }
    \sqrt{t_{\beta}^2-4(-1)^{s_1}} \ &\textrm{ to } \ -\sqrt{t_{\beta}^2-4(-1)^{s_1}}.    
\end{align*}
Since $|\theta_{\alpha, 1}|^{-1}=|\theta_{\alpha, 2}| $, by applying this embedding we get
\[ |\theta_{\alpha, 2}|^{-k_1} |\theta_{\beta, 2}|^{k_2}=1, \]
% i.e.
% \[ \left|\frac{t_{\alpha}-\sqrt{t_{\alpha}^2-4(-1)^{s_0}}}{2}\right|^{-k_1} \left|\frac{t_{\beta}-\sqrt{t_{\beta}^2-4(-1)^{s_1}}}{2}\right|^{k_2}=1. \]
% This is a contradiction as the left-hand side is $<1$ and the right-hand side is $1$.
which is a contradiction as the left hand side is $<1$.
\end{proof}

The following lemma gives an upper bound and lower bound for the terms of a linear recurrence sequence $(G_i)$ in terms of the dominant root of its characteristic polynomial. 
\begin{lemma}\label{G_ub_lb}
    Let $G=(G_{i})_{i \geq 0}$ be a linear recurrence sequence of order $d$ satisfying $G_{i+d}=e_{1}G_{i+d-1}+ \cdots+ e_{d}G_{i}$, $i \geq 0$. If $\gamma_{1}, \ldots, \gamma_{d}$ are the (distinct) roots of the characteristic polynomial with $\gamma_{1}$ as the dominant root, then there exist positive constants $K_{G,1}, K_{G,2}', M_{0}$ depending only on $G$ such that \begin{equation}\label{G_i_ub_lb}
       G_{i}\leq  K_{G,1} \gamma_{1}^{i}\ \textrm{ for }\ i\geq 0  \textrm{ and }\ K_{G,2}'\gamma_{1}^{i} \leq G_{i} \  \textrm{ for } i \geq M_{0}.
    \end{equation} 
\end{lemma}
\begin{proof}
Since $\gamma_{1}, \ldots, \gamma_{d}$ are the distinct roots of the characteristic polynomial $f(x)=x^{d}-e_{1}x^{d-1}-\cdots-e_{d}$, there exist constants $C_{G,1}(>0), \ldots, C_{G,d}$ depending only on $G_{0}, \ldots, G_{d-1}, e_{1}, \ldots, e_{d}$ such that
    \begin{equation}\label{Base_G_poly_rep}
        G_{i} = C_{G,1} \gamma_{1}^{i} - \cdots-C_{G,d} \gamma_{d}^{i},\  i \geq 0
    \end{equation}
 (see \cite[Theorem C.1]{ST86}). Thus
\begin{equation}\nonumber
 G_{i} \leq \gamma_{1}^{i}\left(C_{G,1} + \cdots + C_{G,d}\left|\frac{\gamma_{d}}{\gamma_{1}}\right|^{i} \right).   
\end{equation} \\
Since $\left|\dfrac{\gamma_{j}}{\gamma_1}\right| < 1$ for $2\leq j \leq d$, therefore, 
\begin{equation}\nonumber
 G_{i} \leq K_{G, 1}\gamma_{1}^{i}, \   i \geq 0,  
\end{equation} 
 where $K_{G, 1} = C_{G,1} + \cdots+ |C_{G,d}|$.
Further, there exists $M_{0} \in \mathbb{N}$ depending only on $G$ such that for $i \geq M_{0}$, we have $$ \left|\dfrac{\gamma_{2}}{\gamma_1}\right|^{i} < \frac{C_{G,1}}{2 (|C_{G,2}|+\cdots+|C_{G,d}|)}.$$
Then for $i \geq M_{0}$, we obtain
\begin{align*}
    G_{i} \geq \gamma_{1}^{i} \left| C_{G,1}- \left|\frac{\gamma_{2}}{\gamma_{1}}\right|^{i} \left(|C_{G,2}|+\cdots+|C_{G,d}|\right) \right| &> \frac{C_{G,1}}{2} \gamma_{1}^{i}
    = K_{G,2}'\gamma_{1}^{i},
\end{align*}
where $K_{G,2}'=\nicefrac{C_{G,1}}{2}$.
\end{proof}

The following lemma compares a given positive integer $N$ with a suitable power of $\gamma_{1}$, depending on the highest index in its base-$G$ representation.
\begin{lemma}\label{G-ary_rep_ub_lb}
    Let $G=(G_{i})_{i \geq 0}$ be a linear recurrence sequence of order $d$ satisfying $G_{i+d}=e_{1}G_{i+d-1}+ \cdots+ e_{d}G_{i}$, $i \geq 0$ with $\gamma_{1}$ as its dominant characteristic root. If
    \begin{align*}
             N=f_{1}G_{R_1}+\cdots+f_{s}G_{R_{s}},\ R_1>\cdots>R_{s},
    \end{align*}
is the base-$G$ representation of an integer $N$ with $f_{i} \neq 0$ for $1 \leq i \leq s$, then there exist positive constants  $K_{G,2}, K_{G,3}$ depending only on $G$ such that
\begin{equation}\label{n_G-ary_ub_lb}
 K_{G,2} \gamma_{1}^{R_{1}} \leq N \leq K_{G,3}\gamma_{1}^{R_1}.   
\end{equation}
\end{lemma}
\begin{proof} We have 
\begin{align*}
       N = \sum_{i=1}^{s}f_{i} G_{R_{i}} < G_{R_{1}+1} &< K_{G,1} \gamma_{1}^{R_{1}+1}= K_{G,3} \gamma_{1}^{R_{1}},
  \end{align*} where
  $K_{G,3}=K_{G,1}\gamma_{1}$.
  \noindent
    By Lemma \ref{G_ub_lb}, there exist constants $K_{G,2}'$ and $M_0$, depending only on $G$, such that  for $R_1\geq M_0$, we have  
  \[K_{G,2}' \gamma_{1}^{R_{1}} \leq G_{R_{1}}\leq N=\sum_{i=1}^{s}f_{i} G_{R_{i}}. \]  
  Now, for $R_{1} < M_{0}$, we have 
  $$K_{G,2}'' \gamma_{1}^{R_1} \leq N,$$ where $K_{G, 2}''=\min_{R_{1} < M_{0}} \frac{G_{R_1}}{\gamma_{1}^{R_1}}$.
   Combining, we get 
   $$ K_{G,2} \gamma_{1}^{R_1} \leq G_{R_1} \leq N,$$ where $K_{G,2}=\min(K_{G,2}', K_{G,2}'')$.
\end{proof}

The next lemma provides a sufficient condition for $\theta_{\alpha,1}$ and $\gamma_{1}$ to be multiplicatively independent.
\begin{lemma}
   Let $\alpha$ be a real quadratic irrational number with simple continued fraction expansion
    \[\alpha=[a_{0,0};\ldots,a_{0,r_0-1},\overline{b_{0,0},\ldots,b_{0,s_0-1}}], r_0\geq 0, s_0\geq 1.\] Let $G=(G_{i})_{i \geq 0}$ be a linear recurrence sequence of order $d$  satisfying $G_{i+d}=e_1G_{i+d-1}+\cdots+e_d G_{i}$, $i \geq 0$ with $\gamma_{1}$ as its dominant characteristic root. 
     If $\mathbb{Q}(\sqrt{t_{\alpha}^2-4(-1)^{s_0}}) \neq \mathbb{Q}\left(\gamma_{1}\right)$, then $\theta_{\alpha,1}$ and $\gamma_{1}$ are multiplicatively independent.
\end{lemma}
\begin{proof}
    If $\theta_{\alpha,1}$ and $\gamma_{1}$ are multiplicatively dependent, then there exist integers $k_{1}$ and $k_{2}$ such that $\theta_{\alpha,1}^{k_1} \gamma_{1}^{k_2} =1.$  Since $\theta_{\alpha,1}>1$, $k_{2} \neq 0$. If $k_{1}=0$, then $\gamma_{1}^{k_2}=1$. Recall, from the Introduction (after the definition of a linear recurrence sequence), that $\gamma_1\geq e_1$. Further, by the definition, $e_1>1$ for $d=1$. Hence $\gamma_1>1$ if $d=1$. If $d\geq 2$, then $\gamma_1\neq e_1$ as $f(e_1)=e_{1}^d - e_{1} e_{1}^{d-1}-\cdots - e_{d} <0$, implying that $e_1$ is not a characteristic root. Thus $\gamma_1>e_1\geq 1$ in this case as well.
    Hence $k_{1}\neq 0$. Also, $k_{1}$ and $k_{2}$ cannot both be positive or both be negative. Suppose $k_{1} > 0$ and $k_{2}< 0$.
The number field $ \mathbb{Q}(\gamma_{1}, \sqrt{t_{\alpha}^2-4(-1)^{s_0}} )$ has degree $2d$. There exists an embedding of 
$ \mathbb{Q}(\gamma_{1}, \sqrt{t_{\alpha}^2-4(-1)^{s_0}})$ which takes
\begin{align*}
\gamma_{1} \ &\textrm{ to } \ \gamma_{1}\\
    \noindent \textrm{ and }
    \sqrt{t_{\alpha}^2-4(-1)^{s_0}} \ &\textrm{ to } \ -\sqrt{t_{\alpha}^2-4(-1)^{s_0}}.
\end{align*}
By applying this embedding we get
\[ |\theta_{\alpha, 2}|^{k_1} \left|\frac{1}{\gamma_{1}}\right|^{-k_2}=1, \]
% i.e.
% \[ \left|\frac{t_{\alpha}-\sqrt{t_{\alpha}^2-4(-1)^{s_0}}}{2}\right|^{k_1} \left|\frac{1}{\gamma_{1}}\right|^{-k_2}=1. \]
which is a contradiction as the left hand side is $<1$. The proof is similar if $k_{1} < 0$ and $k_{2}> 0$.

\end{proof}

%=====================================results about height
Next, we recall the definition and some basic properties of the absolute logarithmic height function $h(\cdot)$. If $\delta$ is an algebraic  number with minimal polynomial $f(X)=d_0(X-\delta^{(1)})\cdots(X-\delta^{(d)})\in\mathbb{Z}[X]$, 
    \[
    h(\delta)=\frac{1}{d}\left(\log d_0+\sum_{i=1}^d\max(0,\log|\delta^{(i)}|)\right)
    \]
denotes its absolute logarithmic height.
If $\delta_1,\delta_2$ are algebraic numbers, then
    \begin{align}\label{height_of_product}
 %       h(\alpha\pm\beta)&\leq h(\alpha)+h(\beta)+\log 2,\\
        h(\delta_1\delta_2^{\pm 1})&\leq h(\delta_1)+h(\delta_2)%,\\
 %       h(\alpha^m)&=|m|h(\alpha)\ (m\in\mathbb{Z}).
    \end{align}
(see \cite[Property 3.3]{Wa_00}). The next result gives an upper bound for the absolute logarithmic height of a polynomial in an algebraic number.
\begin{lemma}\cite[Special case of Lemma 3.7]{Wa_00}\label{height_of_poly}
    Let $f(x)\in\mathbb{Z}[x]$ be a non-zero polynomial and let $\delta$ be an algebraic number. Then
    \[
    h(f(\delta))\leq (\deg f)h(\delta)+\log L(f),
    \]
    where $L(f)$ denotes the sum of the absolute values of the coefficients of $f$.
\end{lemma}

The following lemma due to Chim, Pink \& Ziegler \cite{Chim_Ziegler_18} gives a lower bound for the absolute logarithmic height of the product of powers of multiplicatively independent algebraic numbers.
\begin{lemma}\label{CPZ}
    Let $K$ be a number field and suppose that $\delta_1, \delta_2 \in K$ are multiplicatively independent. There exists an effectively computable constant $D_{0} > 0$ such that for $n, m \in \mathbb{Z}$, we have
    $$ h\left(\frac{\delta_{1}^{n}}{\delta_{2}^{m}}\right) \geq D_{0} \max (|n|, |m|).$$
\end{lemma}

%=====================================Baker LFL
We will use the following lower bound for linear forms in logarithms due to Matveev \cite{Mat1, Mat2}. 
\begin{lemma}\label{lfl}
    Let $T$ be a positive integer. Let  $\delta_1,\ldots,\delta_T$ be non-zero real algebraic numbers and $\log\delta_1,\ldots,\log\delta_T$ be some determinations of their complex logarithms. Let $D$ be the degree of the number field generated by $\delta_1,\ldots,\delta_T$ over $\mathbb{Q}$.  For $j=1,\ldots,T$, let $A_j'$ be a real number satisfying
    \[
    A_j'\geq\max\left(Dh(\delta_j),|\log\delta_j|,0.16\right).
    \]
    Let $k_1,\ldots,k_T$ be rational integers. Set
    \[
    B=\max(|k_1|,\ldots,|k_T|) \textrm{ and }\ \Lambda=\delta_1^{k_1}\cdots\delta_T^{k_T}-1.
    \]
    Then
    \[
    \log|\Lambda|>-1.4\times 30^{T+3}T^{4.5}D^{2}\log(eD) A_1'\cdots A_T'\log(eB).
    \]
\end{lemma}
%====================================
We shall need the following lemma due to Peth\"{o} \& de Weger \cite{Petho_Weger_86}.
\begin{lemma}\label{PW}
    Let $a \geq 0$, $c \geq 1$, $g> (\nicefrac{e^2}{c})^{c}$, and let $x \in \mathbb{R}$ be the largest solution of $x= a+g(\log x)^{c}.$ Then, 
    $$ x < 2^{c} (a^{1/c} + g^{1/c} \log (c^{c}g))^{c}.$$
\end{lemma}

\section{Proof of Theorem \ref{TwoOstTech}}\label{Sec_TwoOstTech}
In what follows, $C_1,C_2,\ldots$ denote constants depending only upon\\ $s_0,s_1,a_{0,0},\ldots,a_{0,r_0-1},b_{0,0},\ldots,b_{0,s_0-1},a_{1,0},\ldots,a_{1,r_1-1},b_{1,0},\ldots,b_{1,s_1-1}$.

    Suppose
    \[
n=\epsilon_1q_{\alpha,N_1}+\cdots+\epsilon_kq_{\alpha,N_k}=d_1q_{\beta,M_1}+\cdots+\epsilon_{\ell}q_{\beta,M_{\ell}}
    \]
    where $N_1>\cdots> N_k$ and $M_1>\cdots> M_{\ell}$, are the Ostrowski-$\alpha$ and Ostrowski-$\beta$ representations of $n$, respectively, with $\epsilon_{i} \neq 0$ for $1 \leq i \leq k$ and $d_{i} \neq 0$ for $1 \leq i \leq \ell$. Write
    $N_i=s_0n_i+j_i+r_0$, $0\leq j_i\leq s_0-1$ and
     $M_i=s_1m_i+p_i+r_1$, $0\leq p_i\leq s_1-1$. (Note that  $n_1\geq 2$ as $N_1\geq 2s_0+r_0$.)
    Then
    \begin{equation}\label{rep_n_twobase}
   n= \sum_{i=1}^k \epsilon_{i}q^{(j_i)}_{\alpha,n_i}=\sum_{i=1}^{\ell} d_{i}q^{(p_i)}_{\beta,m_i}.
    \end{equation}
Hence by Lemma \ref{N_ub_lb},% \textcolor{red}{if $n_1\geq N_0(\alpha)$ and $m_1\geq N_0(\beta)$}
    \begin{equation}\label{domterm_samesize}
        c_{\beta, 4}\theta_{\beta,1}^{m_1}\leq n< c_{\alpha, 5}\theta_{\alpha,1}^{n_1}
    \end{equation}
and \begin{equation}\label{domterm_samesize2}
   c_{\alpha, 4}\theta_{\alpha,1}^{n_1}\leq n< c_{\beta,5}\theta_{\beta,1}^{m_1}. 
\end{equation}

\noindent Using \eqref{domterm_samesize}, we get
    \begin{equation}\label{LFL1_B}
    m_1<C_{1} n_1,\ \textrm{ where } \ C_{1} = \dfrac{\log \theta_{\alpha, 1}+\log^+ (\nicefrac{c_{\alpha, 5}}{c_{\beta, 4}})}{\log \theta_{\beta, 1}}
    \end{equation}
    and $\log^+(x)=\max(\log x,0)$.
    
We prove the theorem by induction on $u$. Let $u=4$. Then we have to prove \eqref{tech_min} with $i=4, K_4=L_4=2$. Also, \eqref{tech_max} is trivial. Now, using equations \eqref{binet_rep} and  \eqref{rep_n_twobase}, we have
\[ \epsilon_{1} c_{\alpha,1}^{(j_1)}\theta_{\alpha,1}^{n_1} - \epsilon_{1}c_{\alpha, 2}^{(j_1)}\theta_{\alpha,2}^{n_1}+ \sum_{i=2}^{k} \epsilon_{i}q^{(j_i)}_{\alpha,n_i}= d_{1}c_{\beta,1}^{(p_1)}\theta_{\beta,1}^{m_1} - d_{1}c_{\beta,2}^{(p_1)}\theta_{\beta,2}^{m_1} +\sum_{i=2}^{\ell} d_{i}q^{(p_i)}_{\beta,m_i},\]
i.e.
\begin{equation}\label{diff_eqn_1}
    \epsilon_{1} c_{\alpha,1}^{(j_1)}\theta_{\alpha,1}^{n_1} - d_{1}c_{\beta,1}^{(p_1)}\theta_{\beta,1}^{m_1} =  \epsilon_{1}c_{\alpha, 2}^{(j_1)}\theta_{\alpha,2}^{n_1}  - d_{1}c_{\beta,2}^{(p_1)}\theta_{\beta,2}^{m_1} +\sum_{i=2}^{\ell} d_{i}q^{(p_i)}_{\beta,m_i}- \sum_{i=2}^{k} \epsilon_{i}q^{(j_i)}_{\alpha,n_i}.
\end{equation}
Using equation \eqref{qn_ub_lb}, we have
\begin{align*} 
\left|\sum_{i=2}^{k} \epsilon_{i}q^{(j_i)}_{\alpha,n_i}\right| &\leq \epsilon_{2} c_{\alpha, 3} \theta_{\alpha, 1}^{n_{2}} + \cdots+ \epsilon_{k}   c_{\alpha, 3} \theta_{\alpha, 1}^{n_{k}} \\ &\leq \max(|a_{0,0}|,\ldots, a_{0,r_{0}-1}, b_{0,0},\ldots,b_{0,s_{0}-1}) c_{\alpha, 3} 
 (\theta_{\alpha, 1}^{n_{2}}+ \cdots+1) \\
&\leq C_{2} \theta_{\alpha, 1}^{n_{2}}.
\end{align*}
Similarly,
    \[
    \left|\sum_{i=2}^{\ell} d_{i}q^{(p_i)}_{\beta,m_i}\right|\leq C_{3} \theta_{\beta, 1}^{m_2}.
    \]
From \eqref{diff_eqn_1} and the facts that $|\theta_{\alpha,2}|< 1, |\theta_{\beta,2}|< 1$, we get
        \begin{align*}
           \left|\epsilon_{1} c_{\alpha,1}^{(j_1)}\theta_{\alpha,1}^{n_1}-d_{1}c_{\beta,1}^{(p_1)}\theta_{\beta,1}^{m_1}\right|  &\leq  \left|\epsilon_{1}c_{\alpha, 2}^{(j_1)}\theta_{\alpha,2}^{n_1}-d_{1}c_{\beta,2}^{(p_1)}\theta_{\beta,2}^{m_1}\right|+ \left|\sum_{i=2}^{k}\epsilon_{i}q_{\alpha, n_{i}}^{(j_i)} \right| \\  &+\left|\sum_{i=2}^{\ell}d_{i}q^{(p_i)}_{\beta, m_{i} }\right|\\
    &\leq C_{4} + C_{2} \theta_{\alpha, 1}^{n_2}+ C_{3} \theta_{\beta, 1}^{m_2}. 
        \end{align*}
 
\noindent Thus using \eqref{domterm_samesize2},
    \begin{align}
 \nonumber   \left|\frac{\epsilon_1c_{\alpha,1}^{(j_1)}\theta_{\alpha,1}^{n_1}}{d_{1}c_{\beta,1}^{(p_1)}\theta_{\beta,1}^{m_1}}-1\right|
 &\leq \frac{C_{4} + C_{2} \theta_{\alpha, 1}^{n_2}+ C_{3} \theta_{\beta, 1}^{m_2}}{d_{1}c_{\beta,1}^{(p_1)}\theta_{\beta,1}^{m_1}}\\
 \nonumber
 &\leq \frac{C_{4}}{d_{1}c_{\beta,1}^{(p_1)}} \frac{1}{\theta_{\beta,1}^{m_1}} + \frac{C_{2}}{d_{1}c_{\beta,1}^{(p_1)}} \frac{\theta_{\alpha, 1}^{n_1}}{\theta_{\beta, 1}^{m_1}} \frac{\theta_{\alpha, 1}^{n_2}}{\theta_{\alpha, 1}^{n_1}} + \frac{C_{3}}{d_{1}c_{\beta,1}^{(p_1)}} \frac{\theta_{\beta, 1}^{m_2}}{\theta_{\beta, 1}^{m_1}}\\
 \nonumber
 &\leq \frac{C_{4}}{d_{1}c_{\beta,1}^{(p_1)}} \frac{1}{\theta_{\beta,1}^{m_1}} + \frac{C_{2}}{d_{1}c_{\beta,1}^{(p_1)}} \frac{c_{\beta, 5}}{c_{\alpha,4}} \frac{1}{\theta_{\alpha, 1}^{n_{1}-n_{2}}} + \frac{C_{3}}{d_{1}c_{\beta,1}^{(p_1)}} \frac{1}{\theta_{\beta, 1}^{m_{1}-m_{2}}}\\
 \nonumber
 &\leq C_{5} \max \left\{\frac{1}{\theta^{n_{1}-n_{2}}_{\alpha,1}}, \frac{1}{\theta^{m_{1}-m_{2}}_{\beta,1 }}\right\} \\
 \label{lfl1_rhs}   &\leq \frac{C_{5}}{\min(\theta_{\alpha,1},\theta_{\beta,1})^{\min(n_1-n_2,m_1-m_2)}}.
    \end{align}
Let 
    \[
    \Gamma_{1,1}=\frac{\epsilon_{1}c_{\alpha,1}^{(j_1)}\theta_{\alpha,1}^{n_1}}{d_{1}c_{\beta,1}^{(p_1)}\theta_{\beta,1}^{m_1}}.
    \]
If $\Gamma_{1,1}=1$, then $\dfrac{d_{1} c_{\beta, 1}^{(p_{1})}}{\epsilon_{1} c_{\alpha,1}^{(j_{1})}}=\dfrac{\theta_{\alpha, 1}^{n_{1}}}{\theta_{\beta, 1}^{m_{1}}}$. Taking heights and applying Lemma \ref{CPZ}, we get
$$\max (n_{1}, m_{1}) \leq C_{6}.$$
Therefore
\[
    \min(n_{1} - n_{2},  m_{1}-m_{2}) \leq \max (n_{1}, m_{1})\leq C_{6}\leq 2 C_{6} \log n_{1}.
\]
\noindent
If $\Gamma_{1,1}\neq 1$, then we apply Lemma \ref{lfl} with $T=3$, $D=2$,
    \begin{align*}
        \delta_1&=\theta_{\alpha, 1},\ \delta_2=\theta_{\beta, 1},\ \delta_3=\frac{\epsilon_1c_{\alpha,1}^{(j_1)}}{d_1c_{\beta,1}^{(p_1)}},\\
        k_1&=n_1,\ k_2=-m_1,\ k_3=1,\\
        A_1'&=\log{\theta_{\alpha, 1}},\ A_2'=\log{\theta_{\beta, 1}},\ A_3'=C_{7}.
    \end{align*}
Using \eqref{LFL1_B}, we can take $B=\max(C_{1},1)n_1$. Then 
    \[
  \log|\Gamma_{1,1}-1|>-C_{8}\log(e \max(C_{1},1)n_1).  
    \]
Comparing with \eqref{lfl1_rhs}, we obtain
\[ \min(n_1-n_2,m_1-m_2)< \frac{C_{8}(1+\log \max(C_1, 1)+ \log n_{1})+\log C_{5}}{\log \min(\theta_{\alpha,1}, \theta_{\alpha,2})}.\] Therefore
    \begin{equation*}%\label{lfl1_bound}
    \min(n_1-n_2,m_1-m_2)<C_{9}\log n_1,
    \end{equation*}
where $C_{9}=\dfrac{2 C_{8}\left(2+\log^+ C_1\right)+\log^{+} C_{5}}{\log \min(\theta_{\alpha,1}, \theta_{\alpha,2})}$ since $n_1 \geq 2$.
This proves the base case.\\
Let $u\geq 5$ and suppose there exist integers $K', L', K_4',\ldots, K_{u-1}'$, $L_4',\ldots,$ $ L_{u-1}'$ with
    \begin{align*}
        &2= K_4'\leq\cdots\leq K_{u-1}'=K',\\
        &2= L_4'\leq\cdots\leq L_{u-1}'=L',\\
        &K'+L'=u-1
    \end{align*}
such that \eqref{tech_min} and \eqref{tech_max} hold for $4 \leq i \leq u-1$. We consider three cases.\\
\noindent \textbf{Case 1:} $\boldsymbol{K'>k+1}$\\
Let $(K_i,L_i)=(K'_i,L'_i)$ for $4\leq i \leq u-1$ and let $(K_u,L_u)=(K'+1,L')$. We use the induction hypothesis to prove that the result holds with the above choice.
For $4\leq i \leq u-1$, we have
    \begin{align*}
  \min(n_{1}-n_{K_i},m_{1}-m_{L_i})&=  \min(n_{1}-n_{K'_i},m_{1}-m_{L'_i})\\
  &\leq (C_0\log n_1)^{i-3}\\
\textrm{ and }  \max(n_{1}-n_{K_i-1},m_{1}-m_{L_i-1}) &= \max(n_{1}-n_{K'_i-1},m_{1}-m_{L'_i-1})\\
&\leq (C_0\log n_1)^{i-4}.
    \end{align*}
    Since $n_{j}:= 0$ for $j>k$, we get
    \begin{align*}
  \min(n_{1}-n_{K_u},m_{1}-m_{L_u})&=\min(n_{1},m_{1}-m_{L'})\\
  &=\min(n_{1}-n_{K'_{u-1}},m_{1}-m_{L'_{u-1}})\\
  &\leq (C_0\log n_1)^{(u-1)-3}\leq (C_0\log n_1)^{u-3}.
   \end{align*}
   Further,
   \begin{align*}
  \max(n_{1}-n_{K_{u}-1},m_{1}-m_{L_{u}-1})&=\max(n_{1},m_{1}-m_{L'-1})\\
  &=\max(n_{1}-n_{K'_{u-1}-1},m_{1}-m_{L'_{u-1}-1})\\
  &\leq (C_0\log n_1)^{(u-1)-4}\leq (C_0\log n_1)^{u-4}.
   \end{align*}
% ---------- case 1 ends here
\noindent \textbf{Case 2.}
$\boldsymbol{L'>\ell+1}$\\
Let $(K_i,L_i)=(K'_i,L'_i)$ for $4\leq i \leq u-1$ and let $(K_u,L_u)=(K',L'+1)$.
Already, $m_{j}:=0$ for $j>\ell$. The proof now follows in a manner similar  to Case 1.\\
% ---------- case 2 ends here

\noindent \textbf{Case 3.} $\boldsymbol{K'\leq k+1 \textbf{ and }  L'\leq\ell+1}$\\
By induction hypothesis the equations \eqref{tech_min} \& \eqref{tech_max} hold for $4 \leq i \leq u-1$ and we show that \eqref{tech_min} \& \eqref{tech_max} hold for $i=u$.  Using induction hypothesis for $i=u-1$, we have
\[
\min(n_{1}-n_{K'},m_{1}-m_{L'}) \leq (C_0\log n_1)^{u-4}.
\]
Suppose $\min(n_{1}-n_{K'},m_{1}-m_{L'})=n_{1}-n_{K'}$. Then
\[
    n_{1}-n_{K'} \leq (C_{0}\log n_{1})^{u-4}. 
\]
Let $(K_i,L_i)=(K'_i,L'_i)$ for $4\leq i \leq u-1$ and let $(K_u,L_u)=(K'+1,L')$.
Now, using equations \eqref{binet_rep} and  \eqref{rep_n_twobase}, we have
\begin{align}
 \nonumber \noindent \sum_{i=1}^{K'}(\epsilon_{i} c_{\alpha,1}^{(j_i)}\theta_{\alpha,1}^{n_i} - \epsilon_{i}c_{\alpha, 2}^{(j_i)}\theta_{\alpha,2}^{n_i})&+\sum_{i=K'+1}^{k} \epsilon_{i}q^{(j_i)}_{\alpha,n_i} = \sum_{i=L'}^{\ell} d_{i}q^{(p_i)}_{\beta,m_i}\\
 \nonumber &+\sum_{i=1}^{L'-1}(d_{i}c_{\beta,1}^{(p_i)}\theta_{\beta,1}^{m_i} - d_{i}c_{\beta,2}^{(p_i)}\theta_{\beta,2}^{m_i}),\\
\nonumber i.e.\ \sum_{i=1}^{K'}\epsilon_{i} c_{\alpha,1}^{(j_i)}\theta_{\alpha,1}^{n_i} - \sum_{i=1}^{L'-1}d_{i}c_{\beta,1}^{(p_i)}\theta_{\beta,1}^{m_i}&=\sum_{i=1}^{K'}\epsilon_{i}c_{\alpha, 2}^{(j_i)}\theta_{\alpha,2}^{n_i} - \sum_{i=1}^{L'-1}d_{i}c_{\beta,2}^{(p_i)}\theta_{\beta,2}^{m_i}\\ \label{diff_eqn_2} &+\sum_{i=L'}^{\ell} d_{i}q^{(p_i)}_{\beta,m_i} - \sum_{i=K'+1}^{k} \epsilon_{i}q^{(j_i)}_{\alpha,n_i}.
\end{align}
Using equation \eqref{qn_ub_lb}, we have
% \begin{equation}
%     \nonumber\left|\sum_{i=K'+1}^{k} \epsilon_{i}q^{(j_i)}_{\alpha,n_i}\right| \leq C_{10} \theta_{\alpha, 1}^{n_{K'+1}}
% \end{equation}
% and
%    \begin{equation}
%        \nonumber \left|\sum_{i=L'}^{\ell} d_{i}q^{(p_i)}_{\beta,m_i}\right| \leq C_{11} \theta_{\beta, 1}^{m_{L'}}.
%    \end{equation}  
   \[\left|\sum_{i=K'+1}^{k} \epsilon_{i}q^{(j_i)}_{\alpha,n_i}\right| \leq C_{10} \theta_{\alpha, 1}^{n_{K'+1}} \textrm{ and }\ \left|\sum_{i=L'}^{\ell} d_{i}q^{(p_i)}_{\beta,m_i}\right| \leq C_{11} \theta_{\beta, 1}^{m_{L'}}. \]
Since $|\theta_{\alpha,2}|< 1, |\theta_{\beta,2}|< 1$, using \eqref{diff_eqn_2}, we get
\begin{equation*}
    \left|\sum_{i=1}^{K'}\epsilon_{i} c_{\alpha,1}^{(j_i)}\theta_{\alpha,1}^{n_i} - \sum_{i=1}^{L'-1}d_{i}c_{\beta,1}^{(p_i)}\theta_{\beta,1}^{m_i} \right| \leq C_{12} + C_{10} \theta_{\alpha, 1}^{n_{K'+1}}+C_{11} \theta_{\beta, 1}^{m_{L'}}.
\end{equation*}
Thus using \eqref{domterm_samesize2},
\begin{align}
\nonumber
    \left| \frac{\sum_{i=1}^{K'}\epsilon_{i} c_{\alpha,1}^{(j_i)}\theta_{\alpha,1}^{n_i}}{\sum_{i=1}^{L'-1}d_{i}c_{\beta,1}^{(p_i)}\theta_{\beta,1}^{m_i}}-1\right| &\leq \frac{C_{12} + C_{10} \theta_{\alpha, 1}^{n_{K'+1}}+C_{11} \theta_{\beta, 1}^{m_{L'}}}{\sum_{i=1}^{L'-1}d_{i}c_{\beta,1}^{(p_i)}\theta_{\beta,1}^{m_i}}\\
\nonumber
    &\leq \frac{C_{12} + C_{10} \theta_{\alpha, 1}^{n_{K'+1}}+C_{11} \theta_{\beta, 1}^{m_{L'}}}{d_{1} c_{\beta,1}^{(p_1)} \theta_{\beta,1}^{m_1}}\\
\nonumber
    &\leq \frac{C_{12}}{d_{1} c_{\beta,1}^{(p_1)}\theta_{\beta,1}^{m_{1}}}+ \frac{C_{10}}{d_{1} c_{\beta,1}^{(p_1)}} \frac{\theta_{\alpha,1}^{n_{K'+1}}}{\theta_{\alpha,1}^{n_{1}}} \frac{\theta_{\alpha,1}^{n_{1}}}{\theta_{\beta,1}^{m_1}}+\frac{C_{11}}{d_{1} c_{\beta,1}^{(p_1)}} \frac{\theta_{\beta, 1}^{m_{L'}}}{\theta_{\beta,1}^{m_1}}\\
\nonumber   &\leq C_{13} \frac{1}{\theta_{\beta,1}^{m_{1}}}+ C_{14} \frac{1}{\theta_{\alpha,1}^{n_{1}-n_{K'+1}}} \frac{c_{\beta, 5}}{c_{\alpha,4}}+C_{15} \frac{1}{\theta_{\beta,1}^{m_{1}-m_{L'}}}\\
\nonumber
   & \leq C_{16} \max \left\{\frac{1}{\theta^{n_{1}-n_{K'+1}}_{\alpha,1}}, \frac{1}{\theta^{m_{1}-m_{L'}}_{\beta,1 }}\right\} \\
 \label{lfl1_rhs_1}   &\leq \frac{C_{16}}{\min(\theta_{\alpha,1},\theta_{\beta,1})^{\min(n_1-n_{K'+1},m_{1}-m_{L'})}}.
    \end{align}
Let 
    \[
    \Gamma_{K', L'-1}=\frac{\sum_{i=1}^{K'}\epsilon_{i} c_{\alpha,1}^{(j_i)}\theta_{\alpha,1}^{n_i}}{\sum_{i=1}^{L'-1}d_{i}c_{\beta,1}^{(p_i)}\theta_{\beta,1}^{m_i}}.
    \]
If $\Gamma_{K', L'-1}=1$, then $\frac{\theta_{\alpha,1}^{n_{1}}}{\theta_{\beta,1}^{m_{1}}} = \frac{\sum_{i=1}^{L'-1}d_{i}c_{\beta,1}^{(p_i)}\theta_{\beta,1}^{m_{i}-m_{1}}}{\sum_{i=1}^{K'}\epsilon_{i} c_{\alpha,1}^{(j_i)}\theta_{\alpha,1}^{n_{i}-n_{1}}}$. Applying \eqref{height_of_product}, we get 
\begin{align*}
h\left(\dfrac{\sum_{i=1}^{L'-1}d_{i}c_{\beta,1}^{(p_i)}\theta_{\beta,1}^{m_{i}-m_{1}}}{\sum_{i=1}^{K'}\epsilon_{i} c_{\alpha,1}^{(j_i)}\theta_{\alpha,1}^{n_{i}-n_{1}}}\right) &\leq h\left(\sum_{i=1}^{L'-1}d_{i}c_{\beta,1}^{(p_i)}\theta_{\beta,1}^{m_{i}-m_{1}}\right)\\ &+h\left(\sum_{i=1}^{K'}\epsilon_{i} c_{\alpha,1}^{(j_i)}\theta_{\alpha,1}^{n_{i}-n_{1}}\right).
\end{align*}
We apply Lemma \ref{height_of_poly} with $\delta=\frac{1}{\theta_{\beta, 1}}$, $f=\sum_{i=1}^{L'-1}d_{i}c_{\beta,1}^{(p_i)} x^{m_{1}-m_{i}} $ and $\deg f = m_{1}- m_{L'-1}$ to get
\[ h\left(\sum_{i=1}^{L'-1}d_{i}c_{\beta,1}^{(p_i)}\theta_{\beta,1}^{m_{i}-m_{1}}\right) \leq (m_{1}- m_{L'-1}) \frac{\log \theta_{\beta, 1}}{2} + \log \left(\sum_{i=1}^{L'-1}d_{i}c_{\beta,1}^{(p_i)} \right).\]
Similarly,
\[ h\left(\sum_{i=1}^{K'}\epsilon_{i} c_{\alpha,1}^{(j_i)}\theta_{\alpha,1}^{n_{i}-n_{1}}\right) \leq (n_{1}- n_{K'}) \frac{\log \theta_{\alpha, 1}}{2} + \log \left(\sum_{i=1}^{K'}\epsilon_{i} c_{\alpha,1}^{(j_i)} \right).\]
Now,
\begin{align*}
    \log \left(\sum_{i=1}^{L'-1}d_{i}c_{\beta,1}^{(p_i)} \right) &\leq \log \left(\sum_{i=1}^{\ell}d_{i}c_{\beta,1}^{(p_i)} \right)\\
    &\leq \log \ell + \log \left(\max_{1 \leq i \leq \ell} d_{i} c_{\beta,1}^{(p_i)} \right)
\end{align*}   and
\begin{align*}
     \log \left(\sum_{i=1}^{K'}\epsilon_{i} c_{\alpha,1}^{(j_i)} \right) &\leq \log \left(\sum_{i=1}^{k}\epsilon_{i} c_{\alpha,1}^{(j_i)} \right)\\
    &\leq \log k + \log \left(\max_{1 \leq i \leq k}\epsilon_{i} c_{\alpha,1}^{(j_i)} \right).
\end{align*}
As $k \leq N_{1} \leq  C_{17}\log n_1$ and $\ell \leq M_{1} \leq C_{17}' \log n_1$,
 applying Lemma \ref{CPZ}, there exists an effectively computable constant $D_{0}$, such that  
\begin{align*}
 D_{0} \max(n_{1}, m_{1}) &\leq C_{18}\log n_1 +  (n_{1}-n_{K'}) \frac{\log \theta_{\alpha,1}}{2}+  (m_{1}-m_{L'-1}) \frac{\log \theta_{\beta, 1}}{2}.
\end{align*}
As noted above, 
\begin{equation*}
    n_{1}-n_{K'} \leq (C_{0}\log n_{1})^{u-4}. 
\end{equation*}
\noindent Further, by the induction hypothesis for $i=u-1$, we have
\begin{equation*}
    \max (n_{1}-n_{K'-1}, m_{1}-m_{L'-1}) \leq  (C_{0} \log n_{1})^{u-5}.
\end{equation*}
\noindent Hence 
\begin{equation*}
    \max(n_{1}, m_{1}) \leq  (C_{0} \log n_{1})^{u-4}.
\end{equation*}
Therefore
\begin{align*}
    \min (n_1 -n_{K_u}, m_{1}-m_{L_{u}}) &\leq \max(n_{1}, m_{1}) \\ &\leq  (C_{0} \log n_{1})^{u-4} \leq (C_{0} \log n_{1})^{u-3}
\end{align*}
and
\begin{align*}
     \max (n_1 -n_{K_{u} - 1}, m_{1}-m_{L_{u}-1}) &\leq \max(n_{1}, m_{1}) \leq  (C_{0} \log n_{1})^{u-4}.
\end{align*}
\noindent If $\Gamma_{K', L'-1}\neq 1$, we apply Lemma \ref{lfl} with $T=3$, $D=2$,
\begin{align*}
        \delta_1&=\theta_{\alpha, 1},\ \delta_2=\theta_{\beta, 1},\ \delta_3=\frac{\sum_{i=1}^{K'}\epsilon_{i} c_{\alpha,1}^{(j_i)}\theta_{\alpha,1}^{n_{i}-n_{1}}}{\sum_{i=1}^{L'-1}d_{i}c_{\beta,1}^{(p_i)}\theta_{\beta,1}^{m_{i}-m_{1}}},\\
        k_1&=n_1,\ k_2=-m_1,\ k_3=1,\\
        A_1'&=\log{\theta_{\alpha, 1}},\ A_2'=\log{\theta_{\beta, 1}},\ A_3'=C_{19}(C_{0}\log n_{1})^{u-4}.
    \end{align*}
    % Applying \ref{height_of_product} and properties of absolute logarithmic height on $\delta_{3}$, we get
    % \begin{align*}
    %  h(\delta_{3}) &\leq C_{21} + \log \theta_{\alpha,1}((n_{2}-n_{1})+\cdots+(n_{K'}-n_{1}))+\log \theta_{\beta, 1}((m_{2}-m_{1})+\cdots+(m_{L'-1}-m_{1}))\\
    %  &\leq C_{21} + k \log \theta_{\alpha,1}(n_{K'}-n_{1})+ l \log \theta_{\beta, 1}(m_{L'-1}-m_{1})
    % \end{align*}
    % Using \ref{tech_max} for $i=u-1$ and that $\min(n_{1}-n_{K'},m_{1}-m_{L'})=n_{1}-n_{K'} < (C_{0}\log n_{1})^{u-4}$, we get
    % $$A'_3=C_{22}(C_{0}\log n_{1})^{u-4}$$
Using \eqref{LFL1_B}, we can take $B=\max(C_1,1)n_1$. Then 
    \[
  \log|\Gamma_{K', L'-1}-1|>-C_{20} (C_{0}\log n_{1})^{u-4}\log(e \max(C_{1},1)n_1).  
    \]
Comparing with \eqref{lfl1_rhs_1}, we obtain
    \begin{equation*}%\label{lfl1_bound}
    \min(n_{1}-n_{K'+1},m_{1}-m_{L'})< (C_{0}\log n_1)^{u-3}.
    \end{equation*}
    We will now show that  
    \begin{align*}
        \max(n_{1}-n_{K_{u}-1},m_{1}-m_{L_{u}-1})&=\max(n_{1}-n_{K'},m_{1}-m_{L'-1})\\ &\leq (C_{0} \log n_{1})^{u-4}.
    \end{align*}
Since 
\begin{equation*}
    \min(n_{1}-n_{K'},m_{1}-m_{L'})=n_{1}-n_{K'},
\end{equation*}
we have
\begin{equation*}
    n_{1}-n_{K'} \leq (C_{0}\log n_{1})^{u-4}. 
\end{equation*}
Further, using $i=u-1$, 
$$\max (n_{1}-n_{K'-1}, m_{1}-m_{L'-1}) \leq  (C_{0} \log n_{1})^{u-5}.$$
Thus 
\begin{equation*}
    m_{1}-m_{L'-1} \leq (C_{0}\log n_{1})^{u-5}. 
\end{equation*}
Hence
\begin{equation*}
   \max(n_{1}-n_{K_{u}-1},m_{1}-m_{L_{u}-1}) \leq  (C_{0}\log n_{1})^{u-4}.
\end{equation*}

 % ----------- case3, second subcase starts here

Suppose $\min(n_{1}-n_{K'},m_{1}-m_{L'})=m_{1}-m_{L'}$.\\
Let $(K_i,L_i)=(K'_i,L'_i)$ for $4\leq i \leq u-1$ and let $(K_u,L_u)=(K',L'+1)$.
The proof now follows in a manner similar to the above.

%====================================
\section{Proof of Theorem \ref{OneOstTech}}\label{Sec_OneOstTech}
%====================================
Suppose 
$$n=\epsilon_{1} q_{\alpha, N_{1}}+\cdots+\epsilon_{k} q_{\alpha, N_{k}}=f_{1}G_{R_{1}}+\cdots+f_{s}G_{R_{s}}$$ where $N_{1} > \cdots > N_{k}$ and $R_{1}> \cdots > R_{s}$, are the Ostrowski-$\alpha$ and base-$G$ representations of $n$, respectively, with $\epsilon_{i} \neq 0$ for $1 \leq i \leq k$ and $f_{i} \neq 0$ for $1 \leq i \leq s$. (Note that  $n_1\geq 2$ as $N_1\geq 2s_0+r_0$.)
From \eqref{rep_n_twobase}, $$ n=\sum_{i=1}^{k} \epsilon_{i} q_{\alpha, n_{i}}^{(j_i)} = \sum_{i=1}^{s} f_{i}G_{R_{i}}.$$
Applying \eqref{n_upperbound} and \eqref{n_G-ary_ub_lb}, we get
\begin{equation}\label{BaseG_domterm_samesize}
 K_{G,2} \gamma_{1}^{R_{1}}< n < c_{\alpha, 5} \theta_{\alpha,1}^{n_1}   
\end{equation}
and 
\begin{equation}\label{BaseG_domterm_samesize_2}
 c_{\alpha,4}\theta_{\alpha,1}^{n_1} < n < K_{G,3}\gamma_{1}^{R_1}.  
\end{equation}
Using \eqref{BaseG_domterm_samesize}, we get
    \begin{equation}\label{LFL1_B_1}
    R_1<C_{21} n_1,\ \textrm{ where } \ C_{21} = \dfrac{\log \theta_{\alpha, 1}+\log^+ (\nicefrac{c_{\alpha, 5}}{K_{G, 2}})}{\log \gamma_1}
    \end{equation}
    and $\log^+(x)=\max(\log x,0)$.
    
We prove the theorem by induction on $u$. (The proof is similar to the proof of Theorem \ref{TwoOstTech}. We prove the base case in detail and only outline the rest of the proof.) Let $u=4$. Then we have to prove \eqref{tech_min_2} with $i=4, K_4=L_4=2$. Also, \eqref{tech_max_2} is trivial. Now, using equations \eqref{binet_rep} and  \eqref{Base_G_poly_rep}, we have
\begin{equation}\label{diff_eqn_3}
    \epsilon_{1} c_{\alpha,1}^{(j_1)}\theta_{\alpha,1}^{n_1} - f_{1}C_{G,1}\gamma_{1}^{R_1} =  \epsilon_{1}c_{\alpha, 2}^{(j_1)}\theta_{\alpha,2}^{n_1}- \sum_{i=2}^{d} f_{1} C_{G,i}\gamma_{i}^{R_{1}} +\sum_{i=2}^{s} f_{i} G_{R_{i}} - \sum_{i=2}^{k} \epsilon_{i}q^{(j_i)}_{\alpha,n_i}.
\end{equation}
Using equation \eqref{qn_ub_lb}, we have
\begin{align*} 
\left|\sum_{i=2}^{k} \epsilon_{i}q^{(j_i)}_{\alpha,n_i}\right| 
&\leq C_{22} \theta_{\alpha, 1}^{n_{2}}.
\end{align*}
Similarly, using equation \eqref{G_i_ub_lb}, we have
    \begin{align*}
      \left|\sum_{i=2}^{s} f_{i} G_{R_{i}}\right| &\leq f_{2} K_{G,1}\gamma_{1}^{R_2}+\cdots+f_{s} K_{G,1}\gamma_{1}^{R_s}\\
      &\leq \max(e_{1}, \ldots,e_{d})K_{G,1}(\gamma_{1}^{R_2}+\cdots+1)\\ 
      &\leq C_{23} \gamma_{1}^{R_2}.   
    \end{align*}
   
\noindent Since $|\theta_{\alpha,2}|< 1$ and $|\gamma_{i}|< 1$ for $2 \leq i \leq d$,  using \eqref{diff_eqn_3}, we get
        \begin{align*}
           \left|\epsilon_{1} c_{\alpha,1}^{(j_1)}\theta_{\alpha,1}^{n_1}-f_{1}C_{G,1}\gamma_{1}^{R_1} \right| 
    &\leq C_{24} + C_{22} \theta_{\alpha, 1}^{n_2}+ C_{23}\gamma_{1}^{R_2} . 
        \end{align*}

\noindent Thus using \eqref{BaseG_domterm_samesize_2},
    \begin{align}
 \nonumber   \left|\frac{\epsilon_1c_{\alpha,1}^{(j_1)}\theta_{\alpha,1}^{n_1}}{f_{1}C_{G,1}\gamma_{1}^{R_1}}-1\right|
 &\leq \frac{C_{24} + C_{22} \theta_{\alpha, 1}^{n_2}+ C_{23}\gamma_{1}^{R_2}}{f_{1}C_{G,1}\gamma_{1}^{R_1}}\\
 \nonumber
 &\leq \frac{C_{24}}{f_{1}C_{G,1}} \frac{1}{\gamma_{1}^{R_1}} + \frac{C_{22}}{f_{1}C_{G,1}} \frac{\theta_{\alpha, 1}^{n_1}}{\gamma_{1}^{R_1}} \frac{\theta_{\alpha, 1}^{n_2}}{\theta_{\alpha, 1}^{n_1}} + \frac{C_{23}}{f_{1}C_{G,1}} \frac{\gamma_{1}^{R_2}}{\gamma_{1}^{R_1}}\\
 % \nonumber
 % &\leq \frac{C_{4}}{d_{1}c_{\beta,1}^{(p_1)}} \frac{1}{\theta_{\beta,1}^{m_1}} + \frac{C_{2}}{d_{1}c_{\beta,1}^{(p_1)}} \frac{c_{\beta, 5}}{c_{\alpha,4}} \frac{1}{\theta_{\alpha, 1}^{n_{1}-n_{2}}} + \frac{C_{3}}{d_{1}c_{\beta,1}^{(p_1)}} \frac{1}{\theta_{\beta, 1}^{m_{1}-m_{2}}}\\
 \nonumber
 &\leq C_{25} \max \left\{\frac{1}{\theta^{n_{1}-n_{2}}_{\alpha,1}}, \frac{1}{\gamma_{1}^{R_{1}-R_{2}}}\right\} \\
 \label{lfl1_rhs_2}   &\leq \frac{C_{25}}{\min(\theta_{\alpha,1},\gamma_{1})^{\min(n_1-n_2,R_1-R_2)}}.
    \end{align}
Let 
    \[
    \Gamma_{1,1}'=\frac{\epsilon_{1}c_{\alpha,1}^{(j_1)}\theta_{\alpha,1}^{n_1}}{f_{1}C_{G,1}\gamma_{1}^{R_1}}.
    \]
If $\Gamma_{1,1}'=1$, then $\dfrac{f_{1}C_{G,1}}{\epsilon_{1} c_{\alpha,1}^{(j_{1})}}=\dfrac{\theta_{\alpha, 1}^{n_{1}}}{\gamma_{1}^{R_1}}$. Taking heights and applying Lemma \ref{CPZ}, we get
$$\max (n_{1}, R_{1}) \leq C_{26}.$$
Therefore
\[
    \min(n_{1} - n_{2},  R_{1}-R_{2}) \leq \max (n_{1}, R_{1})\leq C_{26}\leq 2 C_{26} \log n_{1}.
\]
If $\Gamma_{1,1}'\neq 1$, then we apply Lemma \ref{lfl} with $T=3$, $D=\max(2, d)$,
    \begin{align*}
        \delta_1&=\theta_{\alpha, 1},\ \delta_2=\gamma_{1},\ \delta_3=\frac{\epsilon_{1}c_{\alpha,1}^{(j_1)}}{f_{1}C_{G,1}},\\
        k_1&=n_1,\ k_2=-R_1,\ k_3=1,\\
        A_1'&=\log{\theta_{\alpha, 1}},\ A_2'=\log{\gamma_{ 1}},\ A_3'=C_{27}.
    \end{align*}
Using \eqref{LFL1_B_1}, we can take $B=\max(C_{21},1)n_1$. Then 
    \[
  \log|\Gamma_{1,1}'-1|>-C_{28}\log(e \max(C_{21},1)n_1).  
    \]
Comparing with \eqref{lfl1_rhs_2}, we obtain
\[ \min(n_1-n_2,R_1-R_2)< \frac{C_{28}(1+\log \max(C_{21}, 1)+ \log n_{1})+\log C_{25}}{\log \min(\theta_{\alpha,1}, \gamma_{1})}.\] Therefore
    \begin{equation*}%\label{lfl1_bound}
    \min(n_1-n_2,R_1-R_2)<C_{29}\log n_1,
    \end{equation*}
where $C_{29}=\dfrac{2 C_{28}\left(2+\log^+ C_{21}\right)+\log^{+} C_{25}}{\log \min(\theta_{\alpha,1}, \gamma_{1})}$ since $n_1 \geq 2$.
This proves the base case.\\
Let $u\geq 5$ and suppose there exist integers $K', L', K_4',\ldots, K_{u-1}'$, $L_4',\ldots,$ $L_{u-1}'$ with
    \begin{align*}
        &2= K_4'\leq\cdots\leq K_{u-1}'=K',\\
        &2= L_4'\leq\cdots\leq L_{u-1}'=L',\\
        &K'+L'=u-1
    \end{align*}
such that \eqref{tech_min_2} and \eqref{tech_max_2} hold for $4 \leq i \leq u-1$. We consider three cases.\\
\noindent \textbf{Case 1:} $\boldsymbol{K'>k+1}$\\
Let $(K_i,L_i)=(K'_i,L'_i)$ for $4\leq i \leq u-1$ and let $(K_u,L_u)=(K'+1,L')$. The proof follows as in Case 1 of Theorem \ref{TwoOstTech}.\\
% We use the induction hypothesis to prove that the result holds with the above choice.
% For $4\leq i \leq u-1$, we have
%     \begin{align*}
%   \min(n_{1}-n_{K_i},R_{1}-R_{L_i})&=  \min(n_{1}-n_{K'_i},R_{1}-R_{L'_i})\\
%   &\leq (C_{0}'\log n_1)^{i-3}\\
% \textrm{ and }  \max(n_{1}-n_{K_i-1},R_{1}-R_{L_i-1}) &= \max(n_{1}-n_{K'_i-1},R_{1}-R_{L'_i-1})\\
% &\leq (C_{0}'\log n_1)^{i-4}.
%     \end{align*}
%     Since $n_{j}:= 0$ for $j>k$, we get
%     \begin{align*}
%   \min(n_{1}-n_{K_u},R_{1}-R_{L_u})&=\min(n_{1},R_{1}-R_{L'})\\
%   &=\min(n_{1}-n_{K'_{u-1}},R_{1}-R_{L'_{u-1}})\\
%   &\leq (C_{0}'\log n_1)^{(u-1)-3}\leq (C_{0}'\log n_1)^{u-3}.
%    \end{align*}
%    Further,
%    \begin{align*}
%   \max(n_{1}-n_{K_{u}-1},R_{1}-R_{L_{u}-1})&=\max(n_{1},R_{1}-R_{L'-1})\\
%   &=\max(n_{1}-n_{K'_{u-1}-1},R_{1}-R_{L'_{u-1}-1})\\
%   &\leq (C_{0}'\log n_1)^{(u-1)-4}\leq (C_{0}'\log n_1)^{u-4}.
%    \end{align*}
% ---------- case 1 ends here
\noindent \textbf{Case 2.}
$\boldsymbol{L'>s+1}$\\
Let $(K_i,L_i)=(K'_i,L'_i)$ for $4\leq i \leq u-1$ and let $(K_u,L_u)=(K',L'+1)$. The proof follows as in Case 2 of Theorem \ref{TwoOstTech}.\\
% Already, $R_{j}:=0$ for $j>s$. The proof now follows in a manner similar  to Case 1.\\
% ---------- case 2 ends here

\noindent \textbf{Case 3.} $\boldsymbol{K'\leq k+1 \textbf{ and }  L'\leq s+1}$\\
By induction hypothesis the equations \eqref{tech_min_2} \& \eqref{tech_max_2} hold for $4 \leq i \leq u-1$ and we show that \eqref{tech_min_2} \& \eqref{tech_max_2} hold for $i=u$.  Using induction hypothesis for $i=u-1$, we get
\[
\min(n_{1}-n_{K'},m_{1}-m_{L'}) \leq (C_{0}'\log n_1)^{u-4}.
\]
Suppose $\min(n_{1}-n_{K'},m_{1}-m_{L'})=n_{1}-n_{K'}$.\\
Let $(K_i,L_i)=(K'_i,L'_i)$ for $4\leq i \leq u-1$ and let $(K_u,L_u)=(K'+1,L')$.
The rest of the proof can be done in a manner similar to that of the proof of Theorem \ref{TwoOstTech}.

%=======================
\section{Proofs of Theorems \ref{TwoOst} and \ref{OneOst}}\label{Sec_OneTwoOst}
Let $n$ be a positive integer with
    \[
    H_{\alpha}(n)+H_{\beta}(n)\leq M.
    \]
Set $k= H_{\alpha}(n)$ and $\ell=H_{\beta}(n)$. Let
    \[
    n=\epsilon_1q_{\alpha,N_1}+\cdots+\epsilon_kq_{\alpha,N_k},\ N_1>\cdots>N_k,
    \]
be the Ostrowski-$\alpha$ representation of $n$,  with  $\epsilon_{i} \neq 0$ for $1 \leq i \leq k$. Let $n_1=\floor*{\frac{N_1-r_0}{s_0}}$. If $n_{1} \leq 1 $, the result follows trivially. If not, then
$N_1\geq 2s_0+r_0$
and we apply  Theorem\ref{TwoOstTech} with $u=k+\ell+3$ and  \eqref{tech_max} for $i=k+\ell+3$ to get 
% either $K> k+1$ or $L > \ell+1$. We know $n_j =0 $ for $j>k$ and $m_j =0 $ for $j>\ell$,  with $i=k+\ell+3$
\begin{equation*}
n_{1} \leq (C_{0}\log n_{1})^{\ell+k-1}.   
\end{equation*}
We now use Lemma \ref{PW} with $a=0, c=\ell+k-1, g=C_{0}^{\ell+k-1}$. (If necessary, we can replace $C_0$ by $e^2$.) We thus obtain
\begin{align*}
   n_{1} &\leq 2^{M-1} (C_{0} \log ((M-1)^{M-1} C_{0}^{M-1}))^{M-1}\\
   &=  ( 2 C_{0} (M-1) \log ( (M-1) C_{0} ) )^{M-1}\\
   &= ( 2C_{0}\log C_{0} + 2C_{0} \log (M-1))^{M-1} (M-1)^{M-1}.
\end{align*}
Using \eqref{n_upperbound}, we get
\begin{equation*}
    \log n \leq C_{30} n_{1}.
\end{equation*}
Therefore
\begin{align*}
    \log n &\leq C_{30} ( 2C_{0}\log C_{0} + 2C_{0} \log (M-1))^{M-1} (M-1)^{M-1}\\
     &\leq (C M \log M)^{M-1},
\end{align*}
     for $C=4 C_{0} C_{30}^{\nicefrac{1}{M-1}}\log C_{0}$. This completes the proof of Theorem \ref{TwoOst}.\\ Theorem \ref{OneOst} follows similarly from Theorem \ref{OneOstTech}.
\section{Proof of corollary \ref{OneOst_base_b}}\label{sec_6}
\noindent Taking $G_{i} = b^i$ for each $i \geq 0$, we have 
$G_{i+1}= b\, G_{i}$ for all $i \geq 0$. Thus $(G_{i} )_{i \geq 0}$
is a linear recurrence sequence of order $1$ with $\gamma_{1}=b$.
Since $\sqrt{t_{\alpha}^2 - 4(-1)^{s_0}}$ is irrational, we get
$\mathbb{Q}(\sqrt{t_{\alpha}^2 - 4(-1)^{s_0}}) \neq \mathbb{Q}(\gamma_{1})(=\mathbb{Q})$.
 The corollary now follows immediately from Theorem \ref{OneOst}.

 \section*{Acknowledgement}
The authors acknowledge the support of the DST--SERB SRG Grant (SRG/2021/000773). Divyum also thanks BITS Pilani for the support of the OPERA award.
    %=======================
\bibliographystyle{plain}
\bibliography{references}

\begin{thebibliography}{10}

\bibitem{AS}
Jean-Paul Allouche and Jeffrey Shallit.
\newblock {\em Automatic sequences}.
\newblock Cambridge University Press, Cambridge, 2003.
\newblock Theory, applications, generalizations.

\bibitem{AST_022}
Myriam Amri, Lukas Spiegelhofer, and J\"org Thuswaldner.
\newblock R\'epartition jointe dans les classes de r\'esidus de la somme des chiffres pour deux repr\'esentations d'{O}strowski.
\newblock {\em Int. J. Number Theory}, 18(5):955--976, 2022.

\bibitem{dyn}
G.~Barat and P.~Liardet.
\newblock Dynamical systems originated in the {O}strowski alpha-expansion.
\newblock {\em Ann. Univ. Sci. Budapest. Sect. Comput.}, 24:133--184, 2004.

\bibitem{sur}
Val\'{e}rie Berth\'{e}.
\newblock Autour du syst\`eme de num\'{e}ration d'{O}strowski.
\newblock volume 8 (2), pages 209--239. 2001.
\newblock Journ\'{e}es Montoises d'Informatique Th\'{e}orique (Marne-la-Vall\'{e}e, 2000).

\bibitem{BHLS}
Csan{\'a}d Bert{\'o}k, Lajos Hajdu, Florian Luca, and Divyum Sharma.
\newblock {O}n the number of non-zero digits of integers in multi-base representations.
\newblock {\em Publ. Math. Debrecen}, 90(1-2):181--194, 2017.

\bibitem{Be_72}
Jean B\'esineau.
\newblock Ind\'ependance statistique d'ensembles li\'es \`a la fonction ``somme des chiffres''.
\newblock {\em Acta Arith.}, 20:401--416, 1972.

\bibitem{Luca_Bravo_2016}
Jhon~J. Bravo and Florian Luca.
\newblock On the {D}iophantine equation {$F_n+F_m=2^a$}.
\newblock {\em Quaest. Math.}, 39(3):391--400, 2016.

\bibitem{BCM_13}
Yann Bugeaud, Mihai Cipu, and Maurice Mignotte.
\newblock On the representation of {F}ibonacci and {L}ucas numbers in an integer base.
\newblock {\em Ann. Math. Qu\'{e}.}, 37(1):31--43, 2013.

\bibitem{Chim_Ziegler_18}
Kwok~Chi Chim, Istv\'{a}n Pink, and Volker Ziegler.
\newblock On a variant of {P}illai's problem {II}.
\newblock {\em J. Number Theory}, 183:269--290, 2018.

\bibitem{Ziegler_Chim_2018}
Kwok~Chi Chim and Volker Ziegler.
\newblock ``{O}n {D}iophantine equations involving sums of {F}ibonacci numbers and powers of 2''.
\newblock {\em Integers}, 18:Paper No. A99, 30, 2018.

\bibitem{CRT_81}
J.~Coquet, G.~Rhin, and Ph. Toffin.
\newblock Repr\'{e}sentations des entiers naturels et ind\'{e}pendance statistique. {II}.
\newblock {\em Ann. Inst. Fourier (Grenoble)}, 31(1):ix, 1--15, 1981.

\bibitem{EGST}
J.-H. Evertse, K.~Gy\H{o}ry, C.~L. Stewart, and R.~Tijdeman.
\newblock {$S$}-unit equations and their applications.
\newblock In {\em New advances in transcendence theory ({D}urham, 1986)}, pages 110--174. Cambridge Univ. Press, Cambridge, 1988.

\bibitem{Gel_68}
A.~O. Gel$^\prime$fond.
\newblock Sur les nombres qui ont des propri\'et\'es additives et multiplicatives donn\'ees.
\newblock {\em Acta Arith.}, 13:259--265, 1968.

\bibitem{Kim_99}
Dong-Hyun Kim.
\newblock On the joint distribution of $q$-additive functions in residue classes.
\newblock {\em J. Number Theory}, 74(2):307--336, 1999.

\bibitem{LeSh93}
H.~W. Lenstra, Jr. and J.~O. Shallit.
\newblock Continued fractions and linear recurrences.
\newblock {\em Math. Comp.}, 61(203):351--354, 1993.

\bibitem{Lu_00}
Florian Luca.
\newblock Distinct digits in base {$b$} expansions of linear recurrence sequences.
\newblock {\em Quaest. Math.}, 23(4):389--404, 2000.

\bibitem{Mat1}
E.~M. Matveev.
\newblock An explicit lower bound for a homogeneous rational linear form in logarithms of algebraic numbers.
\newblock {\em Izv. Ross. Akad. Nauk Ser. Mat.}, 62(4):81--136, 1998.

\bibitem{Mat2}
E.~M. Matveev.
\newblock An explicit lower bound for a homogeneous rational linear form in logarithms of algebraic numbers. {II}.
\newblock {\em Izv. Ross. Akad. Nauk Ser. Mat.}, 64(6):125--180, 2000.

\bibitem{Mi88}
M.~Mignotte.
\newblock Sur les entiers qui s'\'{e}crivent simplement en diff\'{e}rentes bases.
\newblock {\em European J. Combin.}, 9(4):307--316, 1988.

\bibitem{Os_1922}
Alexander Ostrowski.
\newblock Bemerkungen zur {T}heorie der {D}iophantischen {A}pproximationen.
\newblock {\em Abh. Math. Sem. Univ. Hamburg}, 1(1):77--98, 1922.

\bibitem{Petho_Tichy_89}
A.~Peth\"{o} and R.~F. Tichy.
\newblock On digit expansions with respect to linear recurrences.
\newblock {\em J. Number Theory}, 33(2):243--256, 1989.

\bibitem{Petho_Weger_86}
A.~Peth\"{o} and B.~M.~M. Weger, de.
\newblock Products of prime powers in binary recurrence sequences part i: The hyperbolic case, with an application to the generalized ramanujan-nagell equation.
\newblock {\em Mathematics of Computation}, 47(176):713--727, 1986.

\bibitem{Sc90}
Hans~Peter Schlickewei.
\newblock Linear equations in integers with bounded sum of digits.
\newblock {\em J. Number Theory}, 35(3):335--344, 1990.

\bibitem{SeSt1}
H.~G. Senge and E.~G. Straus.
\newblock {${\rm PV}$}-numbers and sets of multiplicity.
\newblock In {\em Proceedings of the {W}ashington {S}tate {U}niversity {C}onference on {N}umber {T}heory ({W}ashington {S}tate {U}niv., {P}ullman, {W}ash., 1971)}, pages 55--67. Dept. Math., Washington State Univ., Pullman, Wash., 1971.

\bibitem{SeSt2}
H.~G. Senge and E.~G. Straus.
\newblock {${\rm PV}$}-numbers and sets of multiplicity.
\newblock {\em Period. Math. Hungar.}, 3:93--100, 1973.

\bibitem{DS24}
Divyum Sharma.
\newblock On the representation of an imaginary quadratic integer in two different bases.
\newblock {\em Submitted. https://doi.org/10.48550/arXiv.2311.17348}.

\bibitem{DSMa'am_019}
Divyum Sharma.
\newblock Joint distribution in residue classes of the base-{$q$} and {O}strowski digital sums.
\newblock {\em Unif. Distrib. Theory}, 14(2):1--26, 2019.

\bibitem{ST86}
T.~N. Shorey and R.~Tijdeman.
\newblock {\em Exponential Diophantine Equations}.
\newblock Cambridge Tracts in Mathematics. Cambridge University Press, 1986.

\bibitem{Sp23}
L.~Spiegelhofer.
\newblock Collisions of digit sums in bases 2 and 3.
\newblock {\em Isr. J. Math.}, 2023.

\bibitem{St_80}
C.~L. Stewart.
\newblock On the representation of an integer in two different bases.
\newblock {\em J. Reine Angew. Math.}, 319:63--72, 1980.

\bibitem{VZ_21}
Ingrid Vukusic and Volker Ziegler.
\newblock On sums of {F}ibonacci numbers with few binary digits.
\newblock {\em Publ. Math. Debrecen}, 98(1-2):157--181, 2021.

\bibitem{Wa_00}
Michel Waldschmidt.
\newblock {\em Diophantine approximation on linear algebraic groups}, volume 326 of {\em Grundlehren der mathematischen Wissenschaften [Fundamental Principles of Mathematical Sciences]}.
\newblock Springer-Verlag, Berlin, 2000.
\newblock Transcendence properties of the exponential function in several variables.

\bibitem{Ze_72}
E.~Zeckendorf.
\newblock Repr\'esentation des nombres naturels par une somme de nombres de fibonacci ou de nombres de lucas.
\newblock {\em Bull. Soc. Roy. Sci. Li\`ege}, 41:179--182, 1972.

\bibitem{Zi_19}
Volker Ziegler.
\newblock Effective results for linear equations in members of two recurrence sequences.
\newblock {\em Acta Arith.}, 190(2):139--169, 2019.

\end{thebibliography}
\end{document}